\documentclass[a4paper]{article}
\usepackage{amsmath}
\usepackage{amssymb}
\usepackage{amsfonts}
\usepackage{theorem}
\usepackage{blkarray}
\usepackage{arydshln}
\usepackage{multirow}
\usepackage[dvipdfmx]{graphicx}
\makeatletter
\@addtoreset{equation}{section}

\makeatother
\newtheorem{theo}{Theorem}[section]
\newtheorem{defi}{Definition}[section]
\newtheorem{lemm}{Lemma}[section]
\newtheorem{prop}{Proposition}[section]
\newtheorem{exam}{Example}[section]

\newcommand{\Z}{\mathbb Z}
\newcommand{\C}{\mathbb C}
\newcommand{\Q}{\mathbb Q}
\newcommand{\R}{\mathbb R}
\newcommand{\bP}{\mathbb P}
\title{\bf Chow Rings of $\widetilde{Mp}_{0,2}(N,d)$ and $\overline{M}_{0,2}(\bP^{N-1},d)$ and Gromov-Witten Invariants of Projective Hypersurfaces of Degree $1$ and $2$}
\author{Hayato Saito\\
\\
\it Division of Mathematics, Graduate School of Science \\
\it Hokkaido University \\
{\it e-mail address: hayato@math.sci.hokudai.ac.jp}}
\begin{document}

\maketitle


\begin{abstract}
In this paper, we prove formulas that represent two-pointed Gromov-Witten invariant $\langle{\cal O}_{h^{a}}{\cal O}_{h^{b}}\rangle_{0,d}$ of projective hypersurfaces with $d=1,2$ in terms of Chow ring of $\overline{M}_{0,2}(\bP^{N-1},d)$, the moduli spaces of stable maps from 
genus $0$ stable curves to projective space $\bP^{N-1}$. 
Our formulas are based on representation of the intersection number $w({\cal O}_{h^{a}}{\cal O}_{h^{b}})_{0,d}$, which was introduced 
by Jinzenji, in terms of Chow ring of $\widetilde{Mp}_{0,2}(N,d)$, the moduli 
space of quasi maps from $\bP^{1}$ to $\bP^{N-1}$ with two marked points. . 
In order to prove our formulas, we use the results on Chow ring of 
$\overline{M}_{0,2}(\bP^{N-1},d)$, that were derived by A. Musta\c{t}\v{a} and M. Musta\c{t}\v{a}.  
We also present explicit toric data of $\widetilde{Mp}_{0,2}(N,d)$ and 
prove relations of Chow ring of $\widetilde{Mp}_{0,2}(N,d)$.  
\end{abstract}
\section{Introduction}\label{section:intro}

\subsection{Our Aim.}
In computing Gromov-Witten invariants, we usually use \textit{classical mirror symmetry} or \textit{fixed point localization theorem} (\cite{CK} or \cite{Kontsevich}, etc.) especially when we are dealing with basic examples such as projective hypersurfaces.
Furthermore, typical proofs of classical Mirror theorem for toric complete intersections were done by using fixed point localization technique (\cite{CK}, \cite{Givental}, \cite{LLY}). Since localization technique does not need detailed structure of Chow ring of the 
corresponding moduli space, it is still unclear how Gromov-Witten invariants are written in terms of Chow ring of the moduli space.  
In this paper, we prove formulas that represent genus 0 Gromov-Witten invariants of projective hypersurfaces of degree $1$ and $2$ \textit{in terms of Chow ring of the moduli space of stable maps} $\overline{M}_{0,2}(\bP^{N-1},d)$.


\subsection{The Chow Rings of Moduli Spaces of Stable Maps.}
In order to accomplish our program, we need to know detailed structure of the Chow ring of $\overline{M}_{0,2}(\bP^{N-1},d)$.
We mainly refer to Musta\c{t}\v{a}s' results \cite{Mustata1, Mustata2}, and Cox's results \cite{JCox} to obtain information 
we need.

In \cite{JCox}, Cox computed Chow ring of $\overline{M}_{0,2}(\bP^1,2)$.
Its structure is described by using natural basis $D_0,D_1,D_2,H_1,H_2,\psi_1,\psi_2$.
$H_1$ and $H_2$ are pullbacks of hyperplane class with respect to evaluation maps ${\rm ev}_1,\; {\rm ev}_2:\overline{M}_{0,2}(\bP^1,2)\rightarrow \bP^1$.
$\psi_1$ and $\psi_2$ are so called $\psi$-classes of universal curve ${\cal C}\rightarrow \overline{\cal M}_{0,2}(\bP^1,2)$.
$D_i$'s are classes that correspond to loci which parametrize stable maps from nodal curves.
The stable maps that belong to each $D_{i}$ are represented by the following graphs.

\begin{figure}[htbp]
 \begin{center}
  \begin{tabular}{c}

   \begin{minipage}{0.33\hsize}
    \begin{center}
     \includegraphics[width=2.2cm]{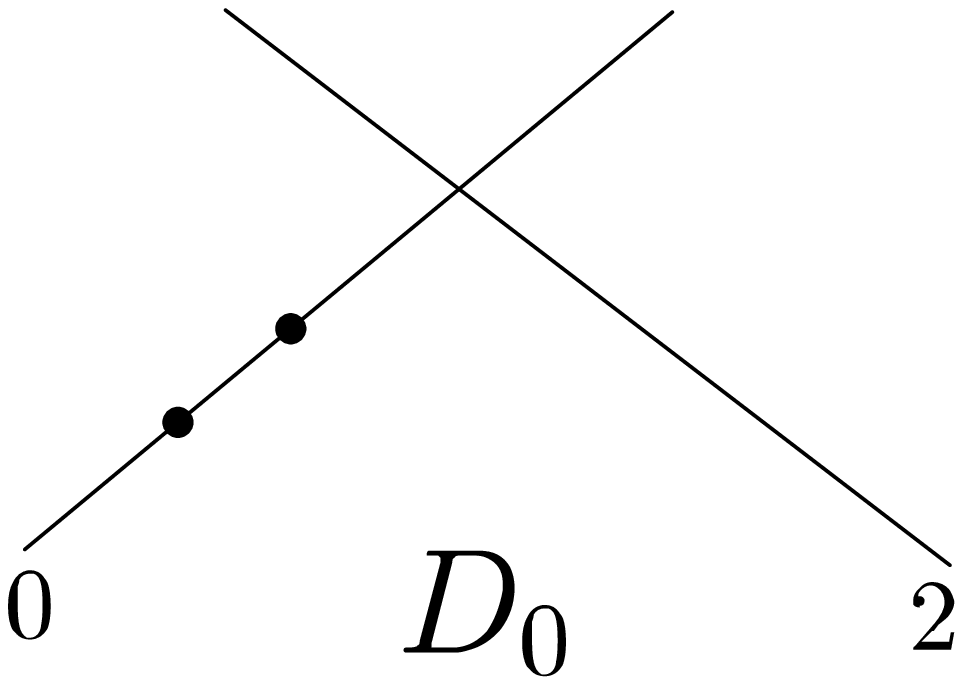}
    \end{center}
   \end{minipage}

   \begin{minipage}{0.33\hsize}
    \begin{center}
     \includegraphics[width=2.2cm]{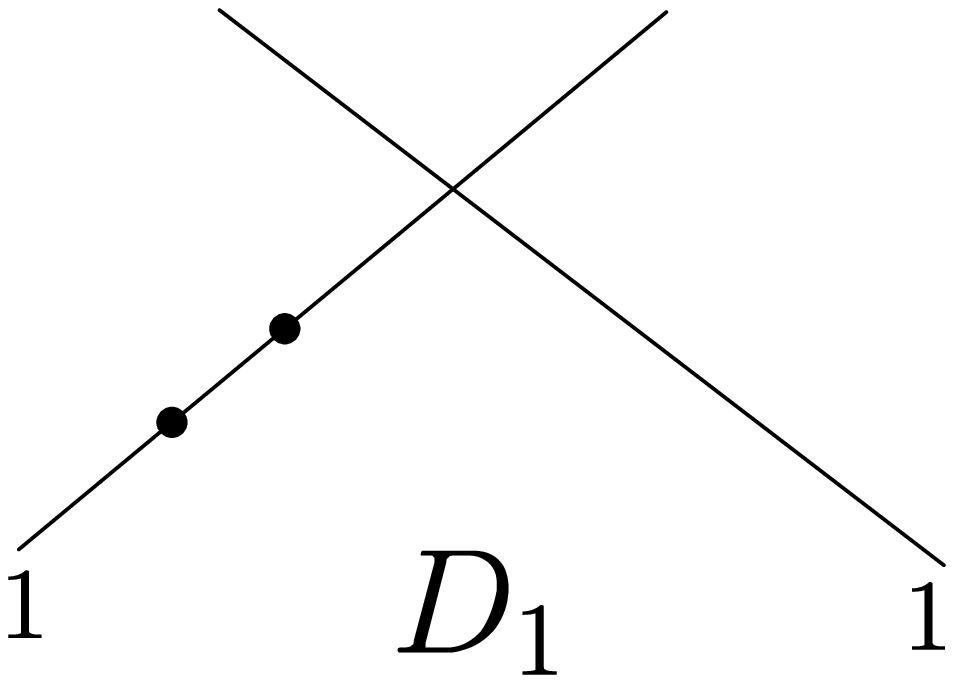}
    \end{center}
   \end{minipage}

   \begin{minipage}{0.33\hsize}
    \begin{center}
     \includegraphics[width=2.2cm]{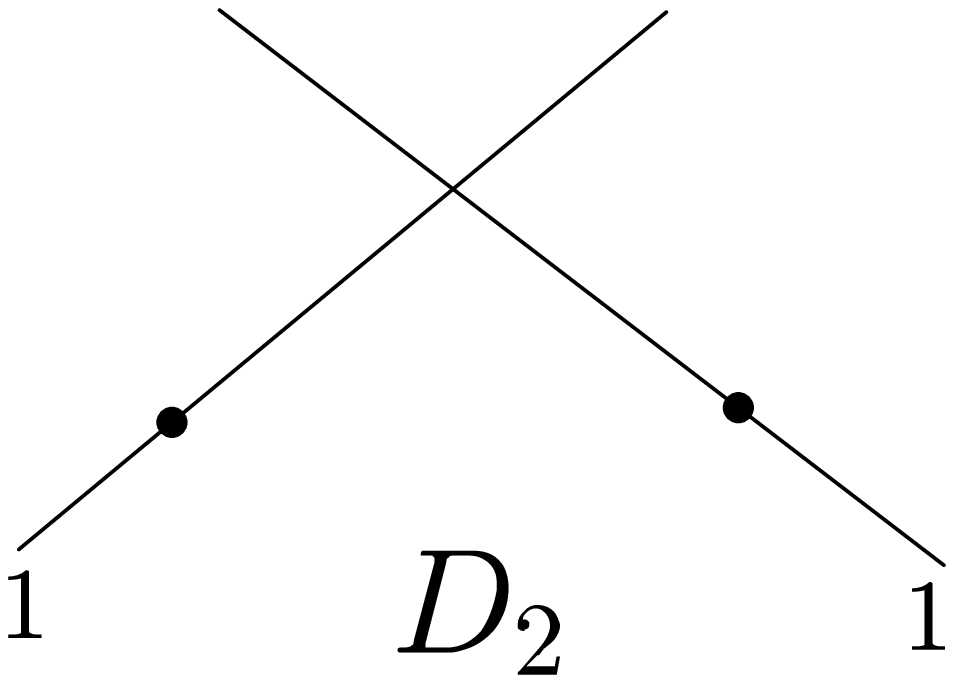}
    \end{center}
   \end{minipage}

  \end{tabular}
 \end{center}
\end{figure}

In \cite{Mustata1}, Musta\c{t}\v{a} and Musta\c{t}\v{a} determined $A^*(\overline{M}_{0,1}(\bP^{N-1},d))$ in general case. They constructed intermediate moduli spaces $\overline{M}_{0,1}(\bP^{N-1},d,k)$ and their substrata $\overline{M}_I^k$, and 
computed extended Chow rings $B^*(\overline{M}_{0,1}(\bP^{N-1},d))$.
$\overline{M}_I^k$ is defined for integer $k\; (0\leq k\leq d)$ and \textit{nested set} $I\subset {\cal P}\backslash \{\emptyset, D\}$, and it parametrizes \textit{$k$-stable maps of $I$-split type} (where ${\cal P}$ is a power set of $D=\{1,2,\dots,d\}$).
The extended Chow ring $B^*(\overline{M}_{0,1}(\bP^{N-1},d))$ is generated by classes associated with $\overline{M}_I^k$'s, and $A^*(\overline{M}_{0,1}(\bP^{N-1},d))$ is given as a subring of $B^*(\overline{M}_{0,1}(\bP^{N-1},d))$ that are invariant under action of symmetric group $S_d$.
In \cite{Mustata2}, they extended their strategy to compute $A^*(\overline{M}_{0,m}(\bP^{N-1},d))$.
For example, generators of $B^*(\overline{M}_{0,2}(\bP^{N-1},2))$ are given by $H$, $\psi$, $T_{\{1_D\}}$, $T_{\{2_D\}}$, $T_{\{1_D,2_D\}}$, $T_{\{1_D,2_M\}}$ and $T_{\{2_D,2_M\}}$.
$H$ is a pullback of hyperplane class via evaluation map ${\rm ev}_1$. $\psi$ coincides with the $\psi_1$ used by Cox.
The others correspond to the following graphs:

\begin{figure}[htbp]
 \begin{center}
  \begin{tabular}{c}

   \begin{minipage}{0.33\hsize}
    \begin{center}
     \includegraphics[width=5cm]{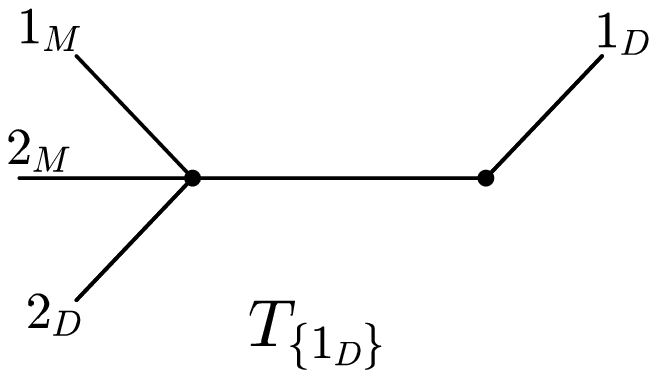}
    \end{center}
   \end{minipage}

   \begin{minipage}{0.33\hsize}
    \begin{center}
     \includegraphics[width=5cm]{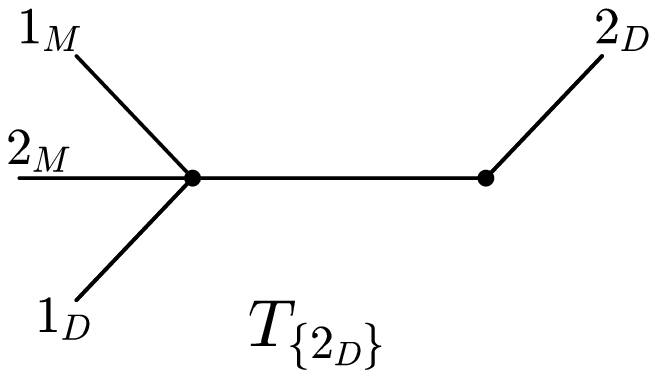}
    \end{center}
   \end{minipage}

   \begin{minipage}{0.33\hsize}
    \begin{center}
     \includegraphics[width=5cm]{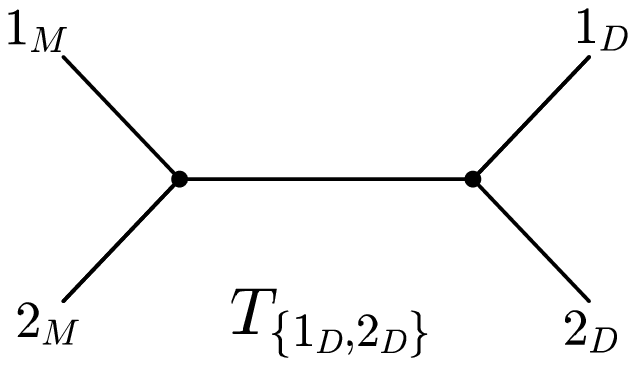}
    \end{center}
   \end{minipage}

  \end{tabular}
 \end{center}
\end{figure}

\begin{figure}[htbp]
 \begin{center}
  \begin{tabular}{c}

   \begin{minipage}{0.33\hsize}
    \begin{center}
     \includegraphics[width=5cm]{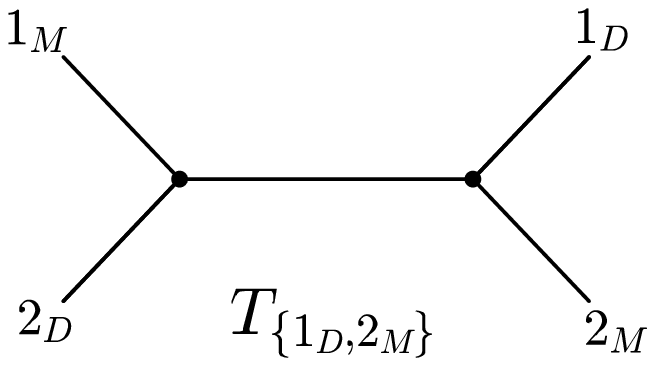}
    \end{center}
   \end{minipage}

   \begin{minipage}{0.33\hsize}
    \begin{center}
     \includegraphics[width=5cm]{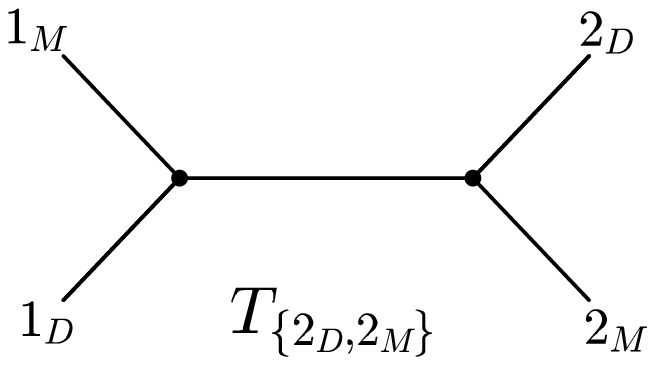}
    \end{center}
   \end{minipage}

  \end{tabular}
 \end{center}
\end{figure}

$T_{\{1_D\}}$ and $T_{\{2_D\}}$ correspond to $D_1$ which was used by Cox.
Similarly, $T_{\{1_D,2_D\}}$ corresponds to $D_0$, $T_{\{1_D,2_M\}}$ and $T_{\{2_D,2_M\}}$ correspond to $D_2$.

Then, basis of $A^*(\overline{M}_{0,2}(\bP^{N-1},2))$ are given by,
\begin{gather*}
 H,\; \psi,\; S_0:=T_{\{1_D,2_D\}},\; S_1:=T_{\{1_D\}}+T_{\{2_D\}},\; S_2:=T_{\{1_D,2_M\}}+T_{\{2_D,2_M\}},\\
 P_1:=T_{\{1_D\}}T_{\{2_D\}},\; P_2:=T_{\{1_D,2_M\}}T_{\{2_D,2_M\}}, \; P_3:=T_{\{2_D\}}T_{\{1_D,2_M\}}+T_{\{1_D\}}T_{\{2_D,2_M\}}.
\end{gather*}


\subsection{Our Motivation 1:  Comparing $\overline{M}_{0,2}(\bP^{N-1},d)$ and $\widetilde{Mp}_{0,2}(N,d)$.}
Our study is influenced by concept of the moduli space $\widetilde{Mp}_{0,2}(N,d)$ of \textit{quasi-maps} from $\bP^{1}$ with 
two marked points to projective space $\bP^{N-1}$, which was introduced by Jinzenji (\cite{Jin1}). This moduli space was also constructed  rigorously by Fontanine and Kim \cite{CK0}. 
In \cite{Jin1}, he presented outline of construction of the moduli space $\widetilde{Mp}_{0,2}(N,d)$, implied that it is a 
toric variety and conjectured generators and relations of its Chow ring.
Although he mentioned the fact that the moduli space $\widetilde{Mp}_{0,2}(N,d)$ is a toric variety, no explicit proof was given.
So, in Section \ref{section:moduli space} of this paper, we prove the following proposition:
\begin{prop}\label{prop:moduli space}
 The space $\widetilde{Mp}_{0,2}(N,d)$ is a toric variety, and its Chow ring is isomorphic to a quotient ring of polynomial ring $\C [H_0,H_1,\dots,H_d]$ modulo an ideal generated by
 \begin{equation}\label{equation:cohomology of moduli space}
  H_0^N,H_1^N(H_0-2H_1+H_2),H_2^N(H_1-2H_2+H_3),\dots,H_{d-1}^N(H_{d-2}-2H_{d-1}+H_d),H_d^N.
 \end{equation}
\end{prop}
We prove it by giving a concrete toric data and using standard theory of toric variety.

From early stage of our study, we have been interested in \textit{bad loci of }$\widetilde{Mp}_{0,2}(N,d)$.
Although $\widetilde{Mp}_{0,2}(N,d)$ parametrizes degree $d$ holomorphic maps from $\bP^1$ to $\bP^{N-1}$ 
with two marked points (see section \ref{section:moduli space} of this paper, or \cite{Jin1}), it has some loci which do \textit{not} correspond to holomorphic maps.
It is given as follows:
let us consider a ``rational map'' $p:\bP^1\rightarrow\bP^{N-1}$ given by
\begin{equation*}
 p([s:t])=[\sum_{j=0}^da_j^0s^{d-j}t^j:\sum_{j=0}^da_j^1s^{d-j}t^j:\cdots:\sum_{j=0}^da_j^Ns^{d-j}t^j].
\end{equation*}
If all the polynomials $\sum_{j=0}^da_j^is^{d-j}t^j$ are divisible by one polynomial $f(s,t)$, then image of $p$
at zero points of $f(s,t)$ cannot be defined.
However, we thought that $\overline{M}_{0,2}(\bP^{N-1},d)$ can be constructed from $\widetilde{Mp}_{0,2}(N,d)$
by successive blow up along these bad loci.
In the $d=1$ case, it is well known that $\overline{M}_{0,2}(\bP^{N-1},1)$ is isomorphic to blow-up of 
$\bP^{N-1}\times \bP^{N-1}$ along its diagonal subset $\Delta$.
Moreover, it can be shown that $\widetilde{Mp}_{0,2}(N,1)$ is isomorphic to $\bP^{N-1}\times \bP^{N-1}$, and its "bad locush is 
givn by the diagonal set $\Delta$.
In this case, $\overline{M}_{0,2}(\bP^{N-1},1)$ \textit{is blow-up of} $\widetilde{Mp}_{0,2}(N,1)$ \textit{along its bad locus}.

Also in general degree $d$, similar claim may be true, but we haven't obtained rigorous proof of this kind of result.
However, we show that Chow ring of $\widetilde{Mp}_{0,2}(N,2)$ is closely related to Chow ring of $\overline{M}_{0,2}(\bP^{N-1},2)$.
To state our result, we transform the basis of $A^*(\overline{M}_{0,2}(\bP^{N-1},2))$ as follows:
\begin{defi}[Key transformation]
 \begin{equation*}
  h_0:=H,\; h_1:=H+\psi,\; h_2:=H+2\psi+S_2.
 \end{equation*}
\end{defi}
Then, we prove the following lemma:
\begin{lemm}\label{lemma:MMp}
 In the ring $A^*(\overline{M}_{0,2}(\bP^{N-1},2))$, the following relations hold. 
 \begin{equation*}
  h_0^N=0,\quad h_1^N(h_0-2h_1+h_2)=0,\quad h_2^N=0.
 \end{equation*}
\end{lemm}
These relations are nothing but the relations of the Chow ring of $\widetilde{Mp}_{0,2}(N,2)$ !

\subsection{Our Motivation 2: Computing GW invariants.}

For degree $k$ projective hypersurface $M_N^k\subset \bP^{N-1}$, Jinzenji introduced intersection number $w({\cal O}_{h^a}{\cal O}_{h^b})_{0,d}$ of $\widetilde{Mp}_{0,2}(N,d)$ defined as follows.
\begin{defi}
 \begin{equation*}
  w({\cal O}_{h^a}{\cal O}_{h^b})_{0,d}:=\int_{\widetilde{Mp}_{0,2}(N,d)} {\rm ev}_1^*(h^a) \wedge {\rm ev}_2^*(h^b) \wedge c_{\rm top}({\cal E}_d^k),
 \end{equation*}
\end{defi}
where ${\cal E}_d^k$ is an orbi-bundle on $\widetilde{Mp}_{0,2}(N,d)$, which corresponds to the condition that the images of quasi maps are contained in degree $k$ hypersurface $M_N^k$.
This orbi-bundle is constructed in \cite{Jin1}.
In \cite{Jin1}, he obtained the result that mirror map of hypersurface $M_N^k$ used in mirror computation of Gromov-Witten invariants is reconstructed as generating function of these intersection numbers, and generalized this framework to the case of toric manifolds with two K\"{a}hler forms.
He also obtained a formula that represents $w({\cal O}_{h^a}{\cal O}_{h^b})_{0,d}$ in terms of Chow ring of $\widetilde{Mp}_{0,2}(N,d)$. 
Let $e^{k}(x,y)=\prod_{j=0}^k(jx+(k-j)y)$. Then $w({\cal O}_{h^a}{\cal O}_{h^b})_{0,d}$ has the following expression.
\begin{equation}\label{equation:w with Chow}
 w({\cal O}_{h^a}{\cal O}_{h^b})_{0,d}=\int_{\widetilde{Mp}_{0,2}(N,d)}H_0^aH_d^b\frac{\prod_{i=0}^{d-1}e^{k}(H_i,H_{i+1})}{\prod_{j=1}^{d-1}(kH_j)}.
\end{equation}
In \cite{Jin3}, he proved mirror formulas that express Gromov-Witten invariant 
$\left<{\cal O}_{h^a} {\cal O}_{h^b}\right>_{0,d}$ of hypersurface $M_N^k$ in terms $w({\cal O}_{h^a}{\cal O}_{h^b})_{0,f}\;(f\leq d)$
and Gromov-Witten invariants of degree lower than $d$ in the $d=1,2,3$ cases. 
For example, in $d=1,2$ cases, GW invariants $\left< {\cal O}_{h^a}{\cal O}_{h^b}\right>_{0,d}$ of degree $k$ hypersurface $M_N^k\subset \bP^{N-1}$ have the following expression.
\begin{align}
 \left<{\cal O}_{h^a} {\cal O}_{h^b}\right>_{0,1}&=w({\cal O}_{h^a}{\cal O}_{h^b})_{0,1}-w({\cal O}_{h^{a+b}}{\cal O}_{h^0})_{0,1}, \label{equation:1 deg GW}\\
({\rm where}&\; a,b\geq 0, a+b=2N-k-3),\nonumber
\end{align}
\begin{align}
 \left<{\cal O}_{h^a} {\cal O}_{h^b}\right>_{0,2}&=w({\cal O}_{h^a}{\cal O}_{h^b})_{0,2}-w({\cal O}_{h^{a+b}}{\cal O}_{h^0})_{0,2}\nonumber \\
&\quad -\frac{1}{k}\left<{\cal O}_{h^a}{\cal O}_{h^b}{\cal O}_{h^{1+k-N}}\right>_{0,1}w({\cal O}_{h^{a+b-N+k}}{\cal O}_{h^0})_{0,1},\label{equation:2 deg GW}\\
({\rm where}&\; a,b\geq 0, a+b=3N-2k-3).\nonumber
\end{align}
We use these formulas to prove our main theorems.

We expected that the equation (\ref{equation:w with Chow}) \textit{can be extended to formulas that express 
the Gromov-Witten invariant $\left<{\cal O}_{h^a} {\cal O}_{h^b}\right>_{0,d}$ in terms of Chow ring of 
$\overline{M}_{0,2}(\bP^{N-1},d)$}, because Chow rings of $\overline{M}_{0,2}(\bP^{N-1},d)$ and $\widetilde{Mp}_{0,2}(N,d)$
are closely related to each other, as can be seen in Lemma \ref{lemma:MMp}.

\subsection{Main Results.}
In Section \ref{section:degree 1} of this paper, we first review the structure of the Chow ring $A^*(\overline{M}_{0,2}(\bP^{N-1},1))$ 
presented in \cite{Mustata2}, and show that three classes  $h_0,h_1,t$ are generators of $A^*(\overline{M}_{0,2}(\bP^{N-1},1))$.
Then we prove the following theorem:
\begin{theo}\label{theorem:Main1}
 For the Gromov-Witten invariant $\left<{\cal O}_{h^a}{\cal O}_{h^b}\right>_{0,1}$ of a hypersurface $M_N^k$, the following 
formula holds:
 \begin{align}
  \int_{\overline{M}_{0,2}(\bP^{N-1},1)}h_0^ah_1^b e^{k}(h_0,h_1+t)=\left<{\cal O}_{h^a} {\cal O}_{h^b}\right>_{0,1}.\nonumber
 \end{align}
\end{theo}

In Section \ref{section: degree 2}, we review the structure of the Chow ring $A^*(\overline{M}_{0,2}(\bP^{N-1},2))$ presented 
in \cite{Mustata2} again, and show that it is generated by five classes $h_0,h_1,h_2,S_0,S_1$ of degree 2 and a class $P_1$ of degree 4.
Furthermore, by using the ring structure in detail, we prove the following:
\begin{theo}\label{theorem:Main2}
 For the Gromov-Witten invariants $\left<{\cal O}_{h^a}{\cal O}_{h^b}\right>_{0,2}$ of hypersurface $M_N^k$, the following 
formula holds:
 \begin{align*}
  &\int_{\overline{M}_{0,2}(\bP^{N-1},2)}(h_0^a-(h_1+S_0)^a)h_2^b \frac{e^{k}(h_0,h_1+S_0)e^{k}(h_1+S_0,h_2+2S_0+S_1)}{k(h_1+S_0)}\\
=&\left<{\cal O}_{h^a}{\cal O}_{h^b}\right>_{0,2}.
 \end{align*}
\end{theo}

Theorem \ref{theorem:Main1} seems to be satisfying to us but Theorem \ref{theorem:Main2} includes somewhat 
strange factor $-(h_1+S_0)^a$.
In fact, $\int_{\overline{M}_{0,2}(\bP^{N-1},2)}h_0^ah_2^b\frac{e^k(h_0,h_1+S_0)e^k(h_1+S_0,h_2+2S_0+S_1)}{k(h_1+S_0)}$ gives us the value  quite close to $\left<{\cal O}_{h^a}{\cal O}_{h^b}\right>_{0,2}$, but does not coincide with it.
In this case, there might be some possibilities to obtain more satisfying formulas. 

\subsection{Organizations of this paper.}
In Section \ref{section:moduli space}, we construct toric data of $\widetilde{Mp}_{0,2}(N,d)$ and compute its Chow ring.
We use some standard theory of toric varieties presented in literatures, e.g. \cite{Ful}, etc.
In Section \ref{section:gwinv}, following \cite{Jin1,Jin2,Jin3}, we derive explicit formulas of Gromov-Witten invariants of projective hypersurface $M_N^k$ in the case of $d=1,2$.
We use these formulas to prove our main results.
In the Section \ref{section:degree 1} and \ref{section: degree 2}, we prove Theorem \ref{theorem:Main1} and \ref{theorem:Main2}.
The section \ref{section: degree 2} contains a proof of Lemma \ref{lemma:MMp}.

\subsection{Notation.}
All varieties considered in this paper are varieties over complex number field.
$M_N^k$ is a hypersurface of degree $k$ in projective space $\bP^{N-1}$.
$M_{m,n}(A)$ is a set of $m\times n$-matrices over ring $A$.

\subsection{Acknowledgements.}
The author would like to thank my supervisor Prof. M. Jinzenji for his support and many helpful discussions.
He is also grateful to Prof. T. Ohmoto for introducing him  Musta\c{t}\v{a}s' works. 
Furthermore, he also would like to thank Prof. B. Kim on comments that led us to our key transformations.  
\section{The Moduli Spase $\widetilde{Mp}_{0,2}(N,d)$}\label{section:moduli space}

In this section we construct toric data of $\widetilde{Mp}_{0,2}(N,d)$ and prove Proposition \ref{prop:moduli space}.
First, we discuss what $\widetilde{Mp}_{0,2}(N,d)$ is (see Section 2.1.1. of \cite{Jin1}).

$\widetilde{Mp}_{0,2}(N,d)$ is a compactified moduli space that parametrizes degree $d$ polynomial maps from 2-pointed $\bP^1$ to $\bP^{N-1}$.
A degree $d$ polynomial map is given by,
\begin{equation*}
 p([s:t])=[{\bf a}_0s^d+{\bf a}_1s^{d-1}t+\cdots+{\bf a}_dt^d],
\end{equation*}
where ${\bf a}_i$'s are vectors in $\C^N$.
Then we introduce uncompactified moduli space $Mp_{0,2}(N,d)$ defined by,
\begin{equation*}
 Mp_{0,2}(N,d):=\{\;({\bf a}_0,{\bf a_1},\dots,{\bf a}_d)\;|\; {\bf a}_0, {\bf a}_{d}\neq {\bf 0}\;\}/(\C^{*})^2.
\end{equation*}
We set the 2-marked points in $\bP^1$ as $0:=[1:0]$ and $\infty:=[0:1]$. Then the condition ${\bf a}_0, {\bf a}_{d}\neq {\bf 0}$
comes from requirement that images of these two marked points are well-defined in $\bP^{N-1}$.
The $(\C^*)^2$ action is induced from the automorphisms of $\bP^1$ which fixes $0$ and $\infty$, and equivalence relation  used in the definition of projective space $\bP^{N-1}$, and it is explicitly written as follows.
\begin{equation*}
 (\mu,\nu) \cdot ({\bf a}_0,{\bf a_1},\dots,{\bf a}_d) = (\mu{\bf a}_0,\mu\nu{\bf a_1},\dots,\mu\nu^{d-1}{\bf a}_{d-1},\mu\nu^d{\bf a}_d).
\end{equation*}
If we use the $(\C^*)^2$ action to turn ${\bf a}_0$ and ${\bf a}_d$ into the points in $\bP^{N-1}$, $Mp_{0,2}(N,d)$ can be describe as follows:
\begin{equation*}
 Mp_{0,2}(N,d)=\{([{\bf a}_0],{\bf a_1},\dots,{\bf a}_{d-1},[{\bf a}_d])\in \bP^{N-1}\times(\C)^{N(d+1)}\times \bP^{N-1}\}/\Z_d.
\end{equation*}
The $\Z_d$-action is given by,
\begin{equation*}
 \zeta_{d} \cdot ([{\bf a}_0],{\bf a_1},\dots,{\bf a}_{d-1},[{\bf a}_d])=([{\bf a}_0],\zeta_d{\bf a_1},\dots,\zeta_d^{d-1}{\bf a}_{d-1},[{\bf a}_d]),
\end{equation*}
where $\zeta_d$ is the $d$-th primitive root of unity.

In order to compactify $Mp_{0,2}(N,d)$, we should add infinite loci corresponding to ${\bf a}_i=\infty$, ($1\leq i\leq d-1$).
In \cite{Jin1}, divisor coordinates $u_i$'s are introduced where the locus ${\bf a}_i=\infty$ is described as zero locus 
of $u_{i}$.
Let $N, d$ be positive integers. The compactified moduli space $\widetilde{Mp}_{0,2}(N,d)$ is defined as follows:
\begin{defi}\label{def:Jin. moduli space}
 \begin{align*}
  &\widetilde{Mp}_{0,2}(N,d) \\
   &:=\{ ({\bf a}_0,\dots,{\bf a}_d,u_1,\dots,u_{d-1}) \in \C^{N(d+1)+d-1} \\
   &\qquad\quad \:|\: {\bf a}_0\neq 0, ({\bf a}_i,u_i)\neq 0\;(1\leq i\leq d-1), {\bf a}_d\neq 0\}/(\C^*)^{d+1},
 \end{align*}
 where the $(\C^*)^{d+1}$-action is given by
 \begin{eqnarray}
  &&(\lambda_0,\dots,\lambda_d)\cdot({\bf a}_0,\dots,{\bf a}_d,u_1,\dots,u_{d-1}) \nonumber \\
 &&=(\lambda_0{\bf a}_0,\lambda_1{\bf a}_1,\dots,\lambda_d{\bf a}_d,\nonumber \\
 &&\qquad\lambda_0^{-1}\lambda_1^2\lambda_2^{-1}u_1,\lambda_1^{-1}\lambda_2^2\lambda_3^{-1}u_2,\dots,\lambda_{d-2}^{-1}\lambda_{d-1}^2\lambda_d^{-1}u_{d-1}). \label{Cd+1-action}
 \end{eqnarray}
\end{defi}


\subsection{Construction of Toric data of $\widetilde{Mp}_{0,2}(N,d)$}
Let
\begin{equation*}
 p_1,p_2,\dots,p_N\in \Z^{N-1}
\end{equation*}
be column vectors which are the 1-skelton of the fan of $\bP^{N-1}$, i.e. 
\begin{equation*}
 (p_1,p_2,\dots,p_{N-1},p_N)=\left(
 \begin{array}{cccccc}
  1        &0         &0        &\cdots  &0         &-1 \\
  0        &1         &0        &\cdots  &0         &-1 \\
  0        &0         &1        &\cdots  &0         &-1\\
  \vdots &\vdots &\vdots &\ddots  &\vdots &\vdots \\
  0        &0         &0        &\cdots  &1         &-1
 \end{array}
 \right)\in M_{N-1,N}(\Z).
\end{equation*}

Next, we introduce $(d+1)$ column vectors
\begin{align*}
 v'_0,v'_1,\dots,v'_d\in \Z^{d-1},
\end{align*}
defined by,
\begin{equation*}
  (v'_0,v'_1,\dots,v'_{d-1},v'_d)=\left(
 \begin{array}{ccccccc}
  -1       &2         &-1      &0         &\cdots  &0          &0 \\
  0        &-1        &2        &-1       & \cdots  &0         &0 \\
  0        &0         &-1       &2         &\cdots  &0         &0\\
  0        &0         &0        &-1        &\cdots  &0         &0\\
  \vdots &\vdots &\vdots &\vdots  &\ddots &\vdots   &0 \\
  0        &0         &0        &0         &\cdots  &-1        &0\\
  0        &0         &0        &0          &\cdots  &2         &-1
 \end{array}
 \right)\in M_{d-1,d+1}(\Z).
\end{equation*}
The $(d-1)\times (d-1)$-submatrix in center of this matrix is Cartan matrix $A_{d-1}$.
Finally, we define column vectors,
\begin{equation*}
 v_{i,j}\;\;(0\leq i\leq d, 1\leq j\leq N), \quad u_k\;\; (1\leq k\leq d-1)
\end{equation*}
as follows:\\
for $j\neq N$,
\begin{equation*}
 v_{i,j}=\begin{blockarray}{*{2}{c}}
  \begin{block}{(c)c}
   {\bf 0}_{N-1}&\\
   \vdots&\\
   p_j& \leftarrow i\\
   \vdots &\\
   {\bf 0}_{N-1}&\\
   {\bf 0}_{d-1}&\\
  \end{block}
 \end{blockarray}
 \in \Z^{(d+1)(N-1)+(d-1)},
\end{equation*}
for $j=N$,
\begin{equation*}
  v_{i,N}=\begin{blockarray}{*{2}{c}}
  \begin{block}{(c)c}
   {\bf 0}_{N-1}&\\
   \vdots&\\
   p_N& \leftarrow i\\
   \vdots &\\
   {\bf 0}_{N-1}&\\
   v'_i&\\
  \end{block}
 \end{blockarray}
 \in \Z^{(d+1)(N-1)+(d-1)}
\end{equation*}
and for $k=1,\dots,d-1$,
\begin{equation*}
 u_k=\begin{blockarray}{*{2}{c}}
  \begin{block}{(c)c}
   {\bf 0}_{N-1}&\\
   \vdots&\\
   {\bf 0}_{N-1}&\\
   -e_k&\\
  \end{block}
 \end{blockarray}
 \in \Z^{(d+1)(N-1)+(d-1)}\end{equation*}
where ${\bf 0}_{N-1}$(resp. ${\bf 0}_{d-1}$) is the zero vector in $\Z^{N-1}$ (resp. $\Z^{d-1}$) and $e_k$ is the $k$-th standard basis of $\Z^{d-1}$.
With this set-up, we define the following toric data:
\begin{defi}\label{definition: fan}
 $\Sigma_{N,d}$ is the fan generated by $\{ v_{i,j} \}_{i=0,\dots,d,\;j=1,\dots,N}$ and $\{ u_k \}_{k=1,\dots, d-1}$.
\end{defi}

\begin{exam}[$N=1$, $d=2$.]
 \begin{equation*}
 v_{0,1}=\begin{blockarray}{*{2}{c}}
  \begin{block}{(c)c}
   1&\\
   0&\\
   0&\\
   0&\\
  \end{block}
 \end{blockarray},
 v_{0,2}=\begin{blockarray}{*{2}{c}}
  \begin{block}{(c)c}
   -1&\\
   0&\\
   0&\\
   -1&\\
  \end{block}
 \end{blockarray},\\
 v_{1,1}=\begin{blockarray}{*{2}{c}}
  \begin{block}{(c)c}
   0&\\
   1&\\
   0&\\
   0&\\
  \end{block}
 \end{blockarray},
 v_{1,2}=\begin{blockarray}{*{2}{c}}
  \begin{block}{(c)c}
   0&\\
   -1&\\
   0&\\
   2&\\
  \end{block}
 \end{blockarray},
\end{equation*}
\begin{equation*}
 v_{2,1}=\begin{blockarray}{*{2}{c}}
  \begin{block}{(c)c}
   0&\\
   0&\\
   1&\\
   0&\\
  \end{block}
 \end{blockarray},
 v_{2,2}=\begin{blockarray}{*{2}{c}}
  \begin{block}{(c)c}
   0&\\
   0&\\
   -1&\\
   -1&\\
  \end{block}
 \end{blockarray},
 u_1=\begin{blockarray}{*{2}{c}}
  \begin{block}{(c)c}
   0&\\
   0&\\
   0&\\
   -1&\\
  \end{block}
 \end{blockarray}.
 \end{equation*}
\end{exam}

Toric variety $X_\Sigma$ defined by the fan $\Sigma=\Sigma_{N,d}$ realizes the variety $\widetilde{Mp}_{0,2}(N,d)$. To see it, we use a standard fact of general theory of toric variety, which describes $X_\Sigma$ as a quotient space of an open subset of $\C^{\Sigma(1)}$. We review here some facts on general fan $\Sigma$ in lattice set $N=\Z^n$ (see Chapter 3 of \cite{CK} or \cite{Ful}).
Let $\Sigma(1)$ be a set of 1-skeltons of $\Sigma$.
Although $\rho\in \Sigma(1)$ is a ray on $N_\R:=N\otimes\R=\R^n$, we identify this ray as generator vector $v_\rho$ of a semi-group $N\cap \rho$.
Let $M:={\rm Hom}_{\Z}(N,\Z)$. We have the following exact sequence for the Chow group $A_{n-1}(X_\Sigma)$:
\begin{equation}\label{exact seq:An-1}
 0\rightarrow M\rightarrow \Z^{\Sigma(1)}\rightarrow A_{n-1}(X_\Sigma)\rightarrow 0.
\end{equation}
where $M\rightarrow \Z^{\Sigma(1)}$ is defined by  $m \mapsto (\left< m\cdot v_\rho\right>)_{\rho\in \Sigma(1)}$ ($\left< m \cdot v_\rho\right>$ is standard inner product),\\ and $\Z^{\Sigma(1)}\rightarrow A_{{\rm dim}(X_{\Sigma})-1}(X_\Sigma)$ is defined by $(m_\rho)\mapsto \sum_{\rho}m_\rho [D_\rho]$ ($[D_\rho]$ is a divisor class associated to $\rho\in \Sigma(1)$).
Next, we define a closed subset $Z(\Sigma)$ of $\C^{\Sigma(1)}$.
The \textit{primitive collection} ${\cal S}$ is a subset of $\Sigma(1)$ which is not the set of 1-dimensional cones of some cone $\sigma\in\Sigma$ but every proper subset of ${\cal S}$ is contained in some cone in the fan. Then, let
\begin{equation*}
 Z(\Sigma):=\bigcup_{{\cal S}{\rm :prim.coll.}} {\bf V}({\cal S}),
\end{equation*}
where ${\bf V}({\cal S})=\{ x\in \C^{\Sigma(1)}\;|\; x_\rho=0, \rho\in {\cal S}\}\subset \C^{\Sigma(1)}$.
Finally, let $G:={\rm Hom}_{\Z}(A_{n-1}(X_\Sigma),\C^*)$.
The group $G$ acts on $\C^{\Sigma(1)}-Z(\Sigma)$ as $g\cdot (x_\rho)_{\rho\in\Sigma(1)}:=(g([D_\rho])x_{\rho})$.
Then the following theorem holds:
\begin{theo}
 If the1-dimensional cones of $\Sigma$ span $N_\R$, then:
 \begin{enumerate}
  \item $X_\Sigma$ is the categorical quotient of $\C^{\Sigma(1)}-Z(\Sigma)$ by $G$,
  \item $X_\Sigma$ is the geometric quotient of $\C^{\Sigma(1)}-Z(\Sigma)$ by $G$ if and only if $X_\Sigma$ is simplicial.
 \end{enumerate}
\end{theo}
See \cite{Cox} for more detail.

We use this theorem in the case of $\Sigma=\Sigma_{N,d}$. Note that $N=M=\Z^{(d+1)(N-1)+(d-1)}$ in this case.   
First, we write down some relations among the vectors in $\Sigma(1)$ to determine the subset $Z(\Sigma)$.
Since we have the relations,
\begin{equation*}
 p_1+p_2+\cdots+p_N={\bf 0}_{N-1},
\end{equation*}
and
\begin{align*}
 v'_0&=-e_1,\\
 v'_1&=-2(-e_1)+(-e_2), \\
 v'_i&=(-e_{i-1})-2(-e_i)+(-e_{i+1})\;\; (2\leq i\leq d-2),\\
 v'_{d-1}&=(-e_{d-2})-2(-e_{d-1}),\\
 v'_d&=(-e_{d-1}),
\end{align*}
we obtain the following relations among the vectors in $\Sigma(1)$:
\begin{align}\label{rel:toric data}
 \Sigma_{j=1}^N v_{0,j}&=u_1,\nonumber \\
 \Sigma_{j=1}^N v_{1,j}&=-2u_1+u_2, \nonumber \\
 \Sigma_{j=1}^N v_{i,j}&=u_{i-1}-2u_i+u_{i+1}\;\; (2\leq i\leq d-2),\\
 \Sigma_{j=1}^N v_{d-1,j}&=u_{d-2}-2u_{d-1}, \nonumber \\
 \Sigma_{j=1}^N v_{d,j}&=u_{d-1}.\nonumber
\end{align}
These relations determine the primitive collections of the fan $\Sigma=\Sigma_{N,d}$ as follows:
\begin{align*}
 &\{ v_{0,1},v_{0,2},\dots, v_{0,N}\},\\
 &\{ v_{1,1},v_{1,2},\dots, v_{1,N},u_1\},\\
 &\{ v_{i,1},v_{i,2},\dots, v_{i,N},u_i\} \;\; (2\leq i\leq d-2), \\
 &\{ v_{d-1,1},v_{d-1,2},\dots, v_{0,N},u_{d-1}\},\\
 &\{ v_{d,1},v_{d,2},\dots, v_{d,N}\}.
\end{align*}
We introduce the following notation, 
\begin{equation*}
 \C^{\Sigma(1)}=\{ x=({\bf a}_0,{\bf a}_1,\dots,{\bf a}_d,u_1,u_2,\dots,u_{d-1}) \; | \; {\bf a}_{i}\in \C^N,\, u_i\in \C \},
\end{equation*}
where ${\bf a}_{i}$ and $u_{i}$ represent $(x_{v_{i,1}},\cdots,x_{v_{i,N}})$ and $x_{u_{i}}$ respectively. 
Then we can easily see that
\begin{equation*}
 Z(\Sigma)=\{ x\in \C^{\Sigma(1)}\;|\; {\bf a}_0=0,\, {\bf a}_d=0,\, ({\bf a}_i,u_i)=0\, (1\leq i\leq d-1)\}.
\end{equation*}
Let us check the $G=(\C^*)^{d+1}$-action on $\C^{\Sigma(1)}$
(note that $|\Sigma(1)|-{\rm rank}(M)=d+1$. then we can see ${\rm rank}(A_{{\rm dim}(X_{\Sigma})-1}(X_\Sigma))=d+1$ from the exact sequence (\ref{exact seq:An-1})). Let $[D_{i,j}]$ (resp. $[U_k]$) be a divisor class that corresponds to 1-skelton $v_{i,j}$ (resp. $u_k$).
By the exact sequence (\ref{exact seq:An-1}) and definition of $v_{i,j}$ and $u_k$, we obtain the following relations on a Chow group $A_{{\rm dim}(X_{\Sigma})-1}(X_{\Sigma})$:
\begin{align*}
 &[D_{i,1}]=[D_{i,2}]=\cdots=[D_{i,N}]\, (0\leq i\leq d),\\
 &[U_k]=-[D_{k-1,N}]+2[D_{k,N}]-[D_{k+1,N}]\, (1\leq k\leq d-1).
\end{align*}
Let $\lambda_i:=g([D_{i,1}])$ ($g\in G$).
The above relation tells us that $g([U_k])=\lambda_{k-1}^{-1}\lambda_k^2\lambda_{k+1}^{-1}$, and the $({\C^*})^{d+1}$-action 
turns out to be,
\begin{align}
 &(\lambda_0,\dots,\lambda_d)\cdot ({\bf a}_0,{\bf a}_1,\dots,{\bf a}_d,u_1,u_2,\dots,u_{d-1}) \nonumber \\
 =&(\lambda_0{\bf a}_0,\lambda_1{\bf a}_1,\dots,\lambda_d{\bf a}_d,\nonumber \\
&\qquad\lambda_0^{-1}\lambda_1^2\lambda_2^{-1}u_1,\lambda_1^{-1}\lambda_2^2\lambda_3^{-1}u_2,\dots,\lambda_{d-2}^{-1}\lambda_{d-1}^2\lambda_d^{-1}u_{d-1}). \nonumber
\end{align}
This action is the same as the $({\C^*})^{d+1}$-action (\ref{Cd+1-action}) in Definition \ref{def:Jin. moduli space}.

\subsection{The Chow ring of $\widetilde{Mp}_{0,2}(N,d)$.}
In this subsection, we prove Proposition \ref{prop:moduli space}.

The Chow ring of $\widetilde{Mp}_{0,2}(N,d)$ is computed by its toric data. To illustrate the recipe of computation, we review general theory of toric varieties (Chap.3 of \cite{CK} or \cite{Ful}, etc.).
Let $\Sigma$ be a simplicial complete fan on a lattice set $N$, and
\begin{equation*}
 \Sigma(1):=\{ \rho_1,\rho_2,\dots, \rho_r\}
\end{equation*}
be an 1-skelton (i.e. a collection of 1-dimensional cones of $\Sigma$).
We identify $\rho_i$ with the generator of semi-group $\rho_i\cap N$ in the same way as the previous subsection, and we denote
it by $v_i$. We define two ideals of $\C[x_1,x_2,\dots,x_r]$.
First one is given by,
\begin{equation*}
 I(\Sigma):=(\sum_{i=1}^r \left<m,v_i\right>\cdot x_i\; |\; m\in M).
\end{equation*}
The other one is the ideal called \textit{Stanley-Reisner ideal} of the corresponding fan $\Sigma$
and is defined by,
\begin{equation*}
 SR(\Sigma):=(x_{i_1}x_{i_2}\cdots x_{i_j}\; |\; \{ v_{i_1},v_{i_2},\dots,v_{i_j}\}\in{\cal S}),
\end{equation*}
where ${\cal S}$ is the primitive collection of the fan $\Sigma$ (see the previous subsection).
Chow ring of a toric variety $X_\Sigma$ is isomorphic to the following quotient ring:
\begin{equation*}
 A^*(X_\Sigma)\cong \C[x_1,x_2,\dots,x_r]/(I(\Sigma)+SR(\Sigma)).
\end{equation*}
This isomorphism enables us to prove Proposition \ref{prop:moduli space}:
\begin{prop}
 $A^*(\widetilde{Mp}_{0,2}(N,d))\cong \C [H_0,H_1,\dots,H_d]/{\cal I}$,\\
where ${\cal I}=(H_0^N,H_d^N,H_k^N(-H_{k-1}+2H_k-H_{k+1}))$
\end{prop}
{\bf proof.}\;\; Let $x_{i,j}$ (resp. $y_k$) be a variable corresponding to a cone $v_{i,j}$ (resp. $u_k$) of the fan $\Sigma_{N,d}$.  Let $\{ e_\alpha\}_{\alpha=1}^{dN+N-2}$ be standard basis of $\Z^{dN+N-2}$. When $1\leq \alpha \leq N-1$, we have,
\begin{align*}
 \left<e_\alpha,v_{i,j}\right>&=
 \begin{cases}
  1\quad (i=0,j=\alpha)\\
  -1\quad (i=0,j=N)\\
  0\quad ({\rm otherwise})
 \end{cases}\\
 \left<e_\alpha,u_k\right>&=0 \; {\rm (for\; all\;}k{\rm )}.
\end{align*}
Hence we obtain $x_{0,\alpha}-x_{0,N}\in I(\Sigma)$, and we can identify $x_{0,\alpha}$ with $x_{0,N}$ in $A^*(\widetilde{Mp}_{0,2}(N,d))$.
In the same manner, when $\ell(N-1)+1 \leq \alpha \leq (\ell+1)(N-1)$ (where $0\leq\ell\leq d$), we have,
\begin{align*}
 \left<e_\alpha,v_{i,j}\right>&=
 \begin{cases}
  1\quad (i=\ell,j=\alpha-\ell(N-1))\\
  -1\quad (i=\ell,j=N)\\
  0\quad ({\rm otherwise})
 \end{cases}\\
 \left<e_\alpha,u_k\right>&=0 \; {\rm (for\; all\;}k{\rm )}.
\end{align*}
Therefore, $x_{\ell,\alpha-\ell(N-1)}=x_{\ell,N}$ in $A^*(\widetilde{Mp}_{0,2}(N,d))$ for $0\leq\ell\leq d$, $\ell(N-1)+1 \leq \alpha \leq (\ell+1)(N-1)$.
Consequently, if we denote $x_{\ell,N}$ by $H_\ell$, we obtain the relation $x_{\ell,j}=H_\ell$ for  $0\leq\ell\leq d$, $1\leq j\leq N-1$.
If $(d+1)(N-1)+1 \leq \alpha \leq (d+1)(N-1)+(d-1)$, we have,
\begin{align*}
 \left<e_\alpha,v_{i,j}\right>&=
 \begin{cases}
  -1\quad (i=\alpha-(d+1)(N-1)-1,j=N)\\
  2\quad (i=\alpha-(d+1)(N-1),j=N)\\
 -1\quad (i=\alpha-(d+1)(N-1)+1,j=N)\\
  0\quad ({\rm otherwise}),
 \end{cases}\\
 \left<e_\alpha,u_k\right>&=
 \begin{cases}
  -1\quad (k=\alpha-(d+1)(N-1))\\
  0\quad ({\rm otherwise}).
 \end{cases}
\end{align*}
Hence  we obtain the relation $-H_{k-1}+2H_k-H_{k+1}-y_k=0$, which yields $y_k=-H_{k-1}+2H_k-H_{k+1}\; (1 \leq k \leq d-1)$.
Accordingly, we have shown that
\begin{equation*}
 \C[x_\rho\;|\; \rho\in\Sigma_{N,d}(1)]/I(\Sigma_{N,d})=\C[H_0,H_1,\dots,H_d].
\end{equation*}
Next, we compute the Stanley-Reisner ideal $SR(\Sigma_{N,d})$.
It is clear that $SR(\Sigma_{N,d})$ coincides with ${\cal I}=(H_0^N,H_d^N,H_k^N(-H_{k-1}+2H_k-H_{k+1}))$ since we have obtained the primitive collection ${\cal S}$ in the previous subsection. $\Box$

\vspace{0.5cm}

In this paper, we use a class:
\begin{equation*}
 d\cdot H_0^{N-1}H_1^NH_2^N\cdots H_{d-1}^NH_d^{N-1}
\end{equation*}
as the volume form of $\widetilde{Mp}_{0,2}(N,d)$, i.e. 
\begin{equation*}
 \int_{\widetilde{Mp}_{0,2}(N,d)}d\cdot H_0^{N-1}H_1^NH_2^N\cdots H_{d-1}^NH_d^{N-1}=1.
\end{equation*}

\section{The Gromov-Witten Invariants of $M_N^k$.}\label{section:gwinv}
In this section, following \cite{Jin1,Jin3}, we derive numerically explicit formulas of the two pointed Gromov-Witten invariants of degree 1 and 2 of $M_N^k$. In the beginning, we define constants which play important roles in the remaining part of this paper.
\begin{defi}
 Let $\ell_i^k$ be the coefficient of $x^{k-i}y^{i+1}$ in $e^k(x,y):=\prod_{j=0}^{k}(jx+(k-j)y)$.
\end{defi}
Then, $e^k(x,y)=\sum_{i=0}^{k-1}\ell_i^kx^{k-i}y^{i+1}$.
Note that $\ell_i^k=\ell_{k-1-i}^k$ and $\ell_i^k=0$ for $i<0$, $k\leq i$.
We describe $\left<{\cal O}_{h^a}{\cal O}_{h^b}\right>_{0,d}$ as a polynomial of $\ell_i^k$ in order to prove our main theorems with  
the aid of the formula (\ref{equation:w with Chow}) and Proposition \ref{prop:moduli space}.

\subsection{The Case of $d=1$.}
Note that $\widetilde{Mp}_{0,2}(N,1)\cong \bP^{N-1}\times \bP^{N-1}$.
We can compute degree $1$ GW invariants of $M_N^k$ with the equation (\ref{equation:1 deg GW}) in Section \ref{section:intro}.
We can compute $w({\cal O}_{h^a}{\cal O}_{h^b})_{0,1}$ by using (\ref{equation:w with Chow}). 
\begin{align*}
 w({\cal O}_{h^a}{\cal O}_{h^b})_{0,1}&=\int H_0^aH_1^b e^k(H_0,H_1)=\int H_0^aH_1^b \sum_{i=0}^{k-1} \ell_i^k H_0^{i+1} H_1^{k-i}\\
 &=\int \sum_{i=0}^{k-1} \ell_i^k H_0^{a+i+1} H_1^{b+k-i}=\int \ell_{N-a-2}^kH_0^{N-1}H_1^{N-1}\\
 &=\ell_{N-a-2}^k.
\end{align*}
Then we obtain,
\begin{align*}
 w({\cal O}_{h^{a+b}}{\cal O}_{h^0})_{0,1}=\ell_{N-a-b-2}^k=\ell_{k-N+1}^k,
\end{align*}
(where we used $a+b=2N-k-3$, and shortened $\int_{\widetilde{Mp}_{0,2}(N,1)}$ to $\int$).
Therefore, we have,
\begin{align}
 \left<{\cal O}_{h^a}{\cal O}_{h^b}\right>_{0,1}&=w({\cal O}_{h^a}{\cal O}_{h^b})_{0,1}-w({\cal O}_{h^{a+b}}{\cal O}_{h^0})_{0,1} \nonumber \\
 &=\ell_{N-a-2}^k-\ell_{k-N+1}^k. \label{GWinv of deg 1}
\end{align}

\subsection{The Case of $d=2$.}
We first present the following formula:
\begin{lemm}
 \begin{equation*}
  \int_{\widetilde{Mp}_{0,2}(N,2)} H_0^\alpha H_1^\beta H_2^\gamma=
  \begin{cases}
   \frac{1}{2^{\beta-N+1}}\binom{\beta-N}{N-\alpha-1} \quad (N \leq \beta \leq 3N-2,\; 0 \leq \alpha,\gamma \leq N-1) \\
   0 \quad ({\rm otherwise}).
  \end{cases}
 \end{equation*}
\end{lemm}
We omit the proof of it because 
it is easily done by using the relation $H_1^N(H_0-2H_1+H_2)=0 \Leftrightarrow H_1^{N+1}=\frac12 H_1^N(H_0+H_2)$ and $H_0^N=H_2^N=0$.

$w({\cal O}_{h^a}{\cal O}_{h^b})_{0,2}$ and $w({\cal O}_{h^{a+b}}{\cal O}_{h^0})_{0,2}$ are computed as follows:
\begin{align*}
 w({\cal O}_{h^a}{\cal O}_{h^b})_{0,2}
 &=\int H_0^aH_2^b \frac{e^k(H_0,H_1)e^k(H_1,H_2)}{kH_1}\\
 &=\int H_0^aH_2^b \frac1k \sum_{i=0}^{k-1}\sum_{j=0}^{k-1} \ell_i^k \ell_j^k H_0^{k-i}H_1^{i+j+1}H_2^{k-j}\\
 &=\int \frac1k \sum_{i,j} \ell_i^k \ell_j^k H_0^{a+k-i}H_1^{i+j+1}H_2^{b+k-j}\\
 &=\frac1k \sum_{i,j} \ell_i^k \ell_j^k \frac{1}{2^{i+j-N+2}}\binom{i+j-N+1}{N-a-k+i-1},\\
 w({\cal O}_{h^{a+b}}{\cal O}_{h^0})_{0,2}
 &=\frac1k \sum_{i,j} \ell_i^k \ell_j^k \frac{1}{2^{i+j-N+2}}\binom{i+j-N+1}{N-(a+b)-k+i-1}\\
 &=\frac1k \sum_{i,j} \ell_i^k \ell_j^k \frac{1}{2^{i+j-N+2}}\binom{i+j-N+1}{-2N+k+i+2}.
\end{align*}

To compute $\left<{\cal O}_{h^a}{\cal O}_{h^b}\right>_{0,2}$ from the equation (\ref{equation:2 deg GW}) in Section \ref{section:intro}, we should evaluate
\\ $w({\cal O}_{h^{a+b-N+k}}{\cal O}_{h^0})_{0,1}$ and $\left<{\cal O}_{h^a}{\cal O}_{h^b}{\cal O}_{h^{1+k-N}}\right>_{0,1}$.
The former is already computed in the previous subsection. Since $a+b=3N-2k-3$ in this case, we obtain,
\begin{equation*}
 w({\cal O}_{h^{a+b-N+k}}{\cal O}_{h^0})_{0,1}=\ell_{k-N+1}^k.
\end{equation*}
The latter can be computed by using Theorem 1 in \cite{Jin2} which computes $n$-\textit{pointed} degree 1 GW invariants of $M_N^k$. If we apply this formula to the case of $n=3$, then
\begin{align*}
 &\left<{\cal O}_{h^a}{\cal O}_{h^b}{\cal O}_{h^{1+k-N}}\right>_{0,1}\\
 &=\int_{\widetilde{Mp}_{0,2}(N,1)} (H_1-H_0) \cdot e^k(H_0,H_1) \cdot H_0^a \frac{H_0^b-H_1^b}{H_0-H_1}\cdot \frac{H_0^{1+k-N}-H_1^{1+k-N}}{H_0-H_1}\\
 &=\int (H_1^{1+k-N}-H_0^{1+k-N})H_0^a\sum_{j=0}^{k-1}\ell_j^k H_0^{k-j}H_1^{j+1} \cdot \sum_{i=0}^{b-1}H_0^iH_1^{b-1-i} \\
 &=\int (H_1^{1+k-N}-H_0^{1+k-N})\sum_{j=0}^{k-1}\sum_{i=0}^{b-1} \ell_j^k H_0^{a+k+i-j}H_1^{b-i+j} \\
 &=\int \sum_{j=0}^{k-1}\sum_{i=0}^{b-1} \ell_j^k H_0^{a+k+i-j}H_1^{k-N+b-i+j+1}
-\int \sum_{j=0}^{k-1}\sum_{i=0}^{b-1} \ell_j^k H_0^{a+2k-N+i-j+1}H_1^{b-i+j}\\
 &=\sum_{i=0}^{b-1}(\ell_{i+a+k-N+1}^k-\ell_{i-b+N-1}^k)\\
 &=\sum_{i=0}^{b-1}(\ell_{i+a+k-N+1}^k-\ell_{i+k-N+1}^k).
\end{align*}
In the last line, we used the symmetry relation $\ell_i=\ell_{k-1-i}$ and inverted the order of terms.
Combining these results, we obtain the following formula.
\begin{align}
 &\left<{\cal O}_{h^a}{\cal O}_{h^b}\right>_{0,2}\nonumber \\
 &=w({\cal O}_{h^a}{\cal O}_{h^b})_{0,2}-w({\cal O}_{h^{a+b}}{\cal O}_{h^0})_{0,2}
 -\frac1k \left<{\cal O}_{h^a}{\cal O}_{h^b}{\cal O}_{h^{1+k-N}}\right>_{0,1}w({\cal O}_{h^{a+b-N+k}}{\cal O}_{h^0})_{0,1}\nonumber \\
 &=\frac1k \sum_{i,j} \ell_i^k \ell_j^k \frac{1}{2^{i+j-N+2}}\Bigg(\binom{i+j-N+1}{N-a-k+i-1}-\binom{i+j-N+1}{-2N+k+i+2}\Bigg)\nonumber \\
&\qquad\qquad -\frac1k \ell_{k-N+1}^k \sum_{i=0}^{b-1}(\ell_{i+a+k-N+1}^k-\ell_{i+k-N+1}^k)\nonumber 
\end{align}
\begin{align}
 &=\frac1k \sum_{i,j} \ell_i^k \ell_j^k \frac{1}{2^{i+j-N+2}}\Bigg(\binom{i+j-N+1}{N-a-k+i-1}-\binom{i+j-N+1}{N-k+j-1}\Bigg)\nonumber \\
 &\qquad\qquad -\frac1k \ell_{k-N+1}^k \sum_{i=0}^{b-1}(\ell_{i+a+k-N+1}^k-\ell_{i+k-N+1}^k). \label{GWinv of deg 2}
\end{align}
\section{Proof of Theorem 1.1.}\label{section:degree 1}
In this section, we will prove Theorem \ref{theorem:Main1} of this paper.
From now on, an integer $n$ is sometimes used as $N-1$:
\begin{equation*}
 n:=N-1.
\end{equation*}
In order to prove the theorem, we use the results on Chow ring of moduli space $\overline{M}_{0,m}(\bP^n,d)$ of stable maps of degree $d$ from genus $0$ stable curve to projective space $\bP^n$.


\subsection{Review of the structure of $A^*(\overline{M}_{0,m}(\bP^n,d))$.}\label{subsection:review of a str of Chow ring}
\begin{theo}[Theorem 5.1. of \cite{Mustata2}]
 Let $M=\{ 1_M, 2_M,\dots, m_M\}$, $D=\{ 1_D,\dots,d_D\}$, $M'=M\backslash\{ 1_M\}$, $D'=M'\sqcup D$ and $d'=|D'|=d+m-1$.

 $B^*(\overline{M}_{0,m}(\bP^n,d))$ is a $\Q-$algebra generated by divisors
 \begin{equation*}
  H\; , \; \psi\; , \; T_h
 \end{equation*}
 for all $h\subset D'$ such that $h\neq \emptyset$ or $\{ i_M\}$ for $i_M\in M'$. Let $T_\emptyset=1$.

 The ideal of relations is generated by:
 \begin{enumerate}
 \renewcommand{\labelenumi}{\rm (\arabic{enumi})}
  \item $H^{n+1}$;
  \item $T_hT_{h'}$ unless $h\cap h'=\emptyset$, or $\emptyset\neq h \subseteq h'$ or $\emptyset\neq h'\subseteq h$;
  \item
   \begin{itemize}
    \item $(m \geq 1)$ $T_hT_{h'}(\psi+\sum_{h\cup h'\subseteq h''}T_{h''})$ for all $h\neq h'$ nonempty;
    \item $(m \geq 2)$ $T_h(\psi+\sum_{h\cup\{ i_M\}\subseteq h'}T_{h'})$ for all $h\neq \emptyset$ and $i_M\in M'\backslash h$;
    \item $(m \geq 3)$ $\psi+\sum_{\{ i_M,j_M\} \subseteq h}T_h$ for all $i_M,j_M\in M'$;
   \end{itemize}
  \item $(m>1)$ $(H+d\psi +\sum_{i_M\in h}|h\cap M'|T_h)^{n+1}$;
  \item $T_h(\sum_{h'\neq h}P(t_{h'})|_{t_{h'}=0}^{t_{h'}=T_{h'}}+\psi^{-1}(H+|^c h_D|\psi)^{n+1})$ for all $h$,
 \end{enumerate}
 \quad where
 \begin{align*}
  &P(t_{h'})=(\psi+\sum_{h''\supset h'}T_{h''}+t_{h''})^{-1}\\
  &[(H+|^c h_D|\psi+\sum_{h''\supset h'}|h_D''\backslash h_D|T_{h''}+|h_D'\backslash h_D|t_{h'})^{n+1}-\\
  &(H+|^c h_D\cap^c h_D'|\psi +\sum_{h''\supset h'}|h_D''\backslash (h_D\cup h_D')|T_{h''})^{n+1}].
 \end{align*}
 Here for any $h\subset D'$, $h_D:=h\cap D$ and $^c h_D:=D\backslash h$.
\label{mustata}
\end{theo}
(We will also use this theorem in the following section, i.e. in the case of degree 2.)
To obtain the Chow ring $A^*(\overline{M}_{0,m}(\bP^n,d))$ from $B^*(\overline{M}_{0,m}(\bP^n,d))$, we have to consider a $S_d\times S_m-$action on $B^*(\overline{M}_{0,m}(\bP^n,d))$.
This action is realized as permutation of $D\cup M$. It is visualized by using the graph that corresponds to $T_h$ in Section \ref{section:intro}.
The Chow ring $A^*(\overline{M}_{0,m}(\bP^n,d))$ is an invariant ring of $B^*(\overline{M}_{0,m}(\bP^n,d))$ under this action (see \cite{Mustata2}).
In the case of $B^*(\overline{M}_{0,2}(\bP^n,1))$,  $A^*(\overline{M}_{0,2}(\bP^n,1))$ coincides with $B^*(\overline{M}_{0,2}(\bP^n,1))$ since this action is trivial. Its generators are
\begin{equation*}
 H,\psi,T_{\{1_D\}}.
\end{equation*}
Let $T:=T_{\{1_D\}}$. Then its relations are the following:
\begin{gather*}
 H^{n+1},\; T\psi,\; (H+\psi)^{n+1},\; \frac{(H+\psi+T)^{n+1}-H^{n+1}}{\psi+T}.
\end{gather*}
Of course, these relations do not contradict the fact that $\overline{M}_{0,2}(\bP^n,1)$ is isomorphic to blow-up of $\bP^n\times\bP^n$ along its diagonal subset $\Delta$. To see it, we define the following transformation:
\begin{defi}[Key transformation]
 \begin{equation*}
  h_0:=H,\; h_1:=H+\psi.
 \end{equation*}
\end{defi}
Then, the above relations change into the following:
\begin{equation*}
 h_0^{n+1},T(h_1-h_0),h_1^{n+1},
\end{equation*}
and
\begin{align}
 \frac{(h_1+T)^{n+1}-h_0^{n+1}}{h_1-h_0+T}.\label{rel5 of 1ptd 1}
\end{align}
At this stage, we introduce a symbol $h$ which satisfies $hT=h_0T=h_1T$ (we have the relation $T(h_1-h_0)$).
Let us expand the last relation:
\begin{align}
 &\frac{(h_1+T)^{n+1}-h_0^{n+1}}{h_1-h_0+T}\nonumber\\
 &=\sum_{i=0}^n h_0^{n-i}(h_1+T)^i =\sum_{i=0}^n h_0^{n-i}\sum_{j=0}^{i}\binom{i}{j}h_1^{i-j}T^j 
 =\sum_{i=0}^n h_0^{n-i}h_1^i +\sum_{i=0}^n\sum_{j=1}^i\binom{i}{j}h^{n-j}T^j \nonumber\\
&=\sum_{i=0}^n h_0^{n-i}h_1^i +\sum_{j=1}^n\sum_{i=j}^n\binom{i}{j}h^{n-j}T^j 
 =\sum_{i=0}^n h_0^{n-i}h_1^i +\sum_{j=1}^n\binom{n+1}{j+1}h^{n-j}T^j.\label{rel5 of 1ptd 2}
\end{align}
In the last line, we used an identity $\sum_{i=j}^n\binom{i}{j}=\binom{n+1}{j+1}$. It can be shown by using $\binom{j}{j}+\binom{j+1}{j}=\binom{j+2}{j+1}$ and Pascal's triangle. 
If we regard $T$ as exceptional divisor of blow-up of $\bP^n\times\bP^n$ along $\Delta$, we can see that these relations coincides  with those of the blow-up (see \cite{Ful2}).
Hence we identify $T$ with the exceptional divisor of the blow-up.

\subsection{Proof of Theorem \ref{theorem:Main1}.}\label{subsection:pf of main1}
In this subsection, we prove the following theorem on degree $1$ Gromov-Witten invariants of a hypersurface $M_{n+1}^k$:
\begin{theo}
 \begin{align*}
  \left<{\cal O}_{h^a}{\cal O}_{h^b}\right>_{0,2}=\int_{\overline{M}_{0,2}(\bP^n,1)}h_0^ah_1^b e^k(h_0,h_1+T),
 \end{align*}
 where $e^k(x,y):=\prod_{j=0}^k (jx+(k-j)y):=\sum_{i=0}^{k-1}\ell_i^kx^{k-i}y^{i+1}$.
\end{theo}
First, we prove the following lemma.
\begin{lemm}\label{lemma1 of 1-ptd}
 \begin{equation*}
  \int_{\overline{M}_{0,2}(\bP^n,1)}h_0^\alpha h_1^\beta (h_1+T)^\gamma=
  \begin{cases}
   1 \; \quad (\alpha = n,\; \gamma \neq n)\\
   -1 \quad (\alpha \neq n,\; \gamma = n)\\
   0 \; \quad ({\rm otherwise})
  \end{cases},
 \end{equation*}
 where $\alpha,\beta,\gamma\geq 0$, $\alpha+\beta+\gamma=2n$.
\end{lemm}
{\bf Proof of lemma \ref{lemma1 of 1-ptd}.}\\
Note that $\int h_0^nh_1^n=1$.
If $\alpha = n$ and $\gamma \neq n$, then $\gamma = n-\beta < n$ and $\beta=n-\gamma>0$. Hence $\alpha+\beta>n$. By the relation $T(h_1-h_0)=0$, the terms that contain $T$ in $h_0^\alpha h_1^\beta (h_1+T)^\gamma$ vanish. Therefore, $h_0^\alpha h_1^\beta (h_1+T)^\gamma=h_0^nh_1^{\beta+\gamma}=h_0^nh_1^n$.
If $\alpha \neq n$ and  $\gamma = n$, then $\beta+\gamma>n$ as above.
Therefore, $h_0^\alpha h_1^\beta (h_1+T)^\gamma=h_0^\alpha h_1^\beta T^n=h^nT^n$.
We obtain the relation $h_0^nh_1^n+h^nT^n=0$ by multiplying the relation (\ref{rel5 of 1ptd 2}) by $h_0^n$, and $h_0^\alpha h_1^\beta (h_1+T)^\gamma=-h_0^nh_1^n$.

Hereinafter, we consider cases when $h_0^\alpha h_1^\beta (h_1+T)^\gamma$ vanishes.
If $\alpha=\gamma=n$ and $\beta=0$, then $h_0^\alpha h_1^\beta (h_1+T)^\gamma=h_0^n(h_1+T)^n=h_0^nh_1^n+h^nT^n=0$.
If $\gamma < n$ and $n-\beta+1 \leq \alpha \leq n-1$, then $\alpha+\beta \geq n+1$. Hence $h_0^\alpha h_1^\beta (h_1+T)^\gamma=h_0^\alpha h_1^{\beta+\gamma}$. However, $\beta+\gamma=2n-\alpha \geq n+1$, and it vanishes.
Finally, we have to consider a case of $\gamma>n$. We multiply the relation (\ref{rel5 of 1ptd 1}) by $h_1-h_0+T$, and obtain the 
relation $(h_1+T)^{n+1}=0$. Therefore $h_0^\alpha h_1^\beta (h_1+T)^\gamma=0$.  $\Box$
\vspace{0.5cm}
\\
\noindent
{\bf Proof of Theorem \ref{theorem:Main1}.}
\begin{align*}
 h_0^ah_1^b e^k(h_0,h_1+T)&=\sum_{i=0}^{k-1}\ell_i^k h_0^{i+1}(h_1+T)^{k-i}=\sum_{i=0}^{k-1}\ell_i^k h_0^{a+i+1}h_1^b(h_1+T)^{k-i}
=(\ell_{n-a-1}^k-\ell_{k-n}^k)h_0^nh_1^n.
\end{align*}
Therefore, $\int_{\overline{M}_{0,2}(\bP^n,1)}h_0^ah_1^b e^k(h_0,h_1+T)=\ell_{n-a-1}^k-\ell_{k-n}^k=\left<{\cal O}_{h^a}{\cal O}_{h^b}\right>_{0,1}$. $\Box$

\section{Proof of Theorem 1.2.}\label{section: degree 2}
In this section, we prove Theorem 1.2 of this paper.
First of all, we have to apply Theorem \ref{mustata} to determine the structure of $A^*(\overline{M}_{0,2}(\bP^n,2))$.
In this case,
\begin{equation*}
 M=\{1_M,2_M\},\;\; D=\{1_D,2_D\},\;\; M'=\{2_M\},\;\; D'=\{1_D,2_D,2_M\}.
\end{equation*}
Therefore, the generators of $B^*(\overline{M}_{0,2}(\bP^n,2))$ are
\begin{equation*}
 H,\;\; \psi,\;\; T_1:=T_{1_D},\;\; T_2:=T_{2_D},\;\; U_1:=T_{1_D,2_M},\;\; U_2:=T_{2_D,2_M},\;\; S_0:=T_{1_D,2_D}.
\end{equation*}
Then, the relations are given as follows:
\begin{align}
 & H^{n+1},\label{rel:1}\\
 & S_0U_1,\;\; S_0U_2,\;\; U_1U_2,\label{rel:2}\\
 & T_1T_2(\psi+S_0)\;\; T_1U_2\psi,\;\; T_2U_1\psi,\label{rel:31}\\
 & T_1(\psi+U_1),\;\; T_2(\psi+U_2),\;\; S_0\psi,\label{rel:32}\\
 & (H+2\psi+U_1+U_2)^{n+1}\label{rel:4},
\end{align}
and the relations (5) in Theorem \ref{mustata}. We take a close look at the relation (5) individually since they have quite 
long expression.

The case of $h=\emptyset$:
\begin{align*}
 &\frac{(H+2\psi+2S_0+U_1+T_1)^{n+1}-(H+\psi+S_0)^{n+1}}{\psi+S_0+U_1+T_1}\\
 -&\frac{(H+2\psi+2S_0+U_1)^{n+1}-(H+\psi+S_0)^{n+1}}{\psi+S_0+U_1}\\
 +&\frac{(H+2\psi+2S_0+U_2+T_2)^{n+1}-(H+\psi+S_0)^{n+1}}{\psi+S_0+U_2+T_2}\\
 -&\frac{(H+2\psi+2S_0+U_2)^{n+1}-(H+\psi+S_0)^{n+1}}{\psi+S_0+U_2}\\
 +&\frac{(H+2\psi+U_1)^{n+1}-(H+\psi)^{n+1}}{\psi+U_1}-\frac{(H+2\psi)^{n+1}-(H+\psi)^{n+1}}{\psi}\\
 +&\frac{(H+2\psi+U_2)^{n+1}-(H+\psi)^{n+1}}{\psi+U_2}-\frac{(H+2\psi)^{n+1}-(H+\psi)^{n+1}}{\psi}\\
 +&\frac{(H+2\psi+2S_0)^{n+1}-H^{n+1}}{\psi+S_0}-\frac{(H+2\psi)^{n+1}-H^{n+1}}{\psi}
 +\frac{(H+2\psi)^{n+1}}{\psi}.
\end{align*}

The cases of $h=\{1_D\}$ and $h=\{2_D\}$:
\begin{align*}
 T_1\Big(& \frac{(H+\psi+S_0+U_2+T_2)^{n+1}-H^{n+1}}{\psi+S_0+U_2+T_2}  & T_2\Big(& \frac{(H+\psi+S_0+U_1+T_1)^{n+1}-H^{n+1}}{\psi+S_0+U_1+T_1}\\
 -&\frac{(H+\psi+S_0+U_2)^{n+1}-H^{n+1}}{\psi+S_0+U_2} & -&\frac{(H+\psi+S_0+U_1)^{n+1}-H^{n+1}}{\psi+S_0+U_1}\\
 +&\frac{(H+\psi+U_2)^{n+1}-H^{n+1}}{\psi+U_2}&+&\frac{(H+\psi+U_1)^{n+1}-H^{n+1}}{\psi+U_1}\\
 -&\frac{(H+\psi)^{n+1}-H^{n+1}}{\psi}&-&\frac{(H+\psi)^{n+1}-H^{n+1}}{\psi}\\
 +&\frac{(H+\psi+S_0)^{n+1}-H^{n+1}}{\psi+S_0}&+&\frac{(H+\psi+S_0)^{n+1}-H^{n+1}}{\psi+S_0}\\
 -&\frac{(H+\psi)^{n+1}-H^{n+1}}{\psi}&-&\frac{(H+\psi)^{n+1}-H^{n+1}}{\psi}\\
 +&\frac{(H+\psi)^{n+1}}{\psi}\Big)&.+&\frac{(H+\psi)^{n+1}}{\psi}\Big).
\end{align*}

The cases of $h=\{1_D,2_M\}$ and $h=\{2_D,2_M\}$:
\begin{align*}
 U_1\Big(&\frac{(H+\psi+S_0+U_2+T_2)^{n+1}-H^{n+1}}{\psi+S_0+U_2+T_2}&U_2\Big(&\frac{(H+\psi+S_0+U_1+T_1)^{n+1}-H^{n+1}}{\psi+S_0+U_1+T_1}\\
 -&\frac{(H+\psi+S_0+U_2)^{n+1}-H^{n+1}}{\psi+S_0+U_2}&-&\frac{(H+\psi+S_0+U_1)^{n+1}-H^{n+1}}{\psi+S_0+U_1}\\
 +&\frac{(H+\psi+U_2)^{n+1}-H^{n+1}}{\psi+U_2}&+&\frac{(H+\psi+U_1)^{n+1}-H^{n+1}}{\psi+U_1}\\
 -&\frac{(H+\psi)^{n+1}-H^{n+1}}{\psi}&-&\frac{(H+\psi)^{n+1}-H^{n+1}}{\psi}\\
 +&\frac{(H+\psi+S_0)^{n+1}-H^{n+1}}{\psi+S_0}&+&\frac{(H+\psi+S_0)^{n+1}-H^{n+1}}{\psi+S_0}\\
 -&\frac{(H+\psi)^{n+1}-H^{n+1}}{\psi}&-&\frac{(H+\psi)^{n+1}-H^{n+1}}{\psi}\\
 +&\frac{(H+\psi)^{n+1}}{\psi}\Big).&+&\frac{(H+\psi)^{n+1}}{\psi}\Big).
\end{align*}

The relation of $h=\{1_D,2_D\}$ is trivial.

We can simplify these relations by using the relations from (1) to (4). We write down the simplest form:
\begin{align}
 &\sum_{i=0}^n(H+\psi+S_0)^{n-i}((H+2\psi+2S_0+U_1+T_1)^i+(H+2\psi+2S_0+U_2+T_2)^i\nonumber \\
 &\qquad\qquad\qquad\qquad-(H+2\psi+2S_0+U_1)^i-(H+2\psi+2S_0+U_2)^i)\nonumber \\
 +&\sum_{i=0}^n(H+\psi)^{n-i}((H+2\psi+S_2)^i-(H+2\psi)^i)\nonumber\\
 +&2\sum_{i=0}^nH^{n-i}(H+2\psi+2S_0)^i,\label{rel:5e}
\end{align}
\begin{align}
 T_1\sum_{i=0}^nH^{n-i}((H+U_2+T_2)^i+(H+\psi+S_0)^i-H^i),\label{rel:511}\\
 T_2\sum_{i=0}^nH^{n-i}((H+U_1+T_1)^i+(H+\psi+S_0)^i-H^i),\label{rel:512}
\end{align}
\begin{align}
 U_1\sum_{i=0}^nH^{n-i}((H+T_2)^i-H^i+(H+\psi)^i),\label{rel:521}\\
 U_2\sum_{i=0}^nH^{n-i}((H+T_1)^i-H^i+(H+\psi)^i).\label{rel:522}
\end{align}

Generators of $A^*(\overline{M}_{0,2}(\bP^n,2))$ are given by,
\begin{align*}
 H,\;\; \psi,\;\; S_0,\;\; S_1:=T_1+T_2,\;\; S_2:=U_1+U_2,\;\; \\
 P_1:=T_1T_2,\;\; P_2:=U_1U_2,\;\; P_3:=T_1U_2+T_2U_1.
\end{align*}
But we don't have to use $P_2$ and $P_3$.
$P_2=U_{1}U_{2}$ vanishes by the relation (\ref{rel:2}). 
On the other hand, by using the relation (\ref{rel:32}), we obtain,
\begin{equation*}
 0=T_1(\psi+U_1)+T_1(\psi+U_1)=S_1\psi+S_1S_2-P_3,
\end{equation*}
Hence $P_3$ is represented in terms of $\psi$, $S_1$ and $S_2$.

We can derive some useful relations from them by applying the following transformation:
\begin{defi}[Key transformation]
 \begin{equation*}
  h_0:=H,\;\; h_1:=H+\psi,\;\; h_2:=H+2\psi+S_2.
 \end{equation*}
\end{defi}
Then the generators of $A^*(\overline{M}_{0,2}(\bP^n,2))$ turn into $h_0,h_1,h_2,S_0,S_1$ and $P_1$.

\begin{lemm}\label{rel:mp}
 \begin{equation*}
  h_0^{n+1}=0,\;\; h_2^{n+1}=0,\;\; h_1^{n+1}(h_0-2h_1+h_2)=0.
 \end{equation*}
\end{lemm}
{\bf Proof.} \; The first two are clear from (\ref{rel:1}) and (\ref{rel:4}). If we multiply by $\psi$ the relation (\ref{rel:521}) and (\ref{rel:522}), and add the resulting relations, then we obtain the third relation as follows:
\begin{align*}
 &\psi \cdot {\rm (\ref{rel:521})} + \psi \cdot {\rm (\ref{rel:522})}&\\
 =&\psi U_1 \sum_{i=0}^nH^{n-i}(H+\psi)^i + \psi U_2 \sum_{i=0}^nH^{n-i}(H+\psi)^i & {\rm (by \; (\ref{rel:31}))}\\
 =&S_2 ((H+\psi)^{n+1}-H^{n+1})&\\
 =&(h_0-2h_1+h_2)h_1^{n+1}.&\Box
\end{align*}
These relations are interesting since they are the same as the ones of $A^*(\widetilde{Mp}_{0,2}(N,2))$. Hence we use
the same volume form as the one of $\widetilde{Mp}_{0,2}(N,2)$. 
\begin{equation*}
 \int_{\overline{M}_{0,2}(\bP^n,2)}2h_0^nh_1^{n+1}h_2^n=1.
\end{equation*}
If there is no room for misunderstanding, we abbreviate $\int_{\overline{M}_{0,2}(\bP^n,2)}$ as $\int$.
We write down the relations from (\ref{rel:1}) to (\ref{rel:4}) in terms of new generators.
\begin{lemm}\label{lemm:14}
 \begin{align*}
  P_3=S_1(h_2-h_1),\;\;P_1(h_2-h_0)=0,\;\; P_1(h_1-h_0+S_0)=0,\\
  S_0(h_1-h_0)=0,\;\; S_0(h_1-h_2)=0,\;\; S_1(h_1-h_0)(h_1-h_2)=0.
 \end{align*}
\end{lemm}
We have to compute the following intersection number in order to prove of Theorem \ref{theorem:Main2}:
\begin{equation}\label{int num}
 \int h_0^\alpha (h_1+S_0)^\beta (h_2+2S_0+S_1)^\gamma h_2^\delta.
\end{equation}
For this purpose, we set
\begin{equation*}
 g_0:=h_0,\;\; g_1:=h_1+S_0,\;\; g_2:=h_2+2S_0+S_1,
\end{equation*}
and derive relations among $g_0$, $g_1$, and $g_2$.

\begin{lemm}\label{lemm:ggg}
 \begin{equation*}
  g_0^{n+1}=0,\;\;g_1^{n+1}(g_0-2g_1+g_2)=0,\;\; (g_1-g_0)\sum_{i=0}^n(g_0^i+g_2^i)g_1^{n-i}=0,\;\; (g_1-g_0)g_2^{n+1}=0.
 \end{equation*}
\end{lemm}
{\bf Proof.} \; The first relation is trivial.
The second relation is equivalent to
\begin{equation*}
 S_1(h_1+S_0)^{n+1}=0
\end{equation*}
by lemma \ref{rel:mp}. It is shown by adding (\ref{rel:511}) and (\ref{rel:512}), and by multiplying the resulting expression by $(\psi+S_0)$.
\begin{align*}
 0=&(\psi+S_0)((\ref{rel:511}) + (\ref{rel:512}))\\
 =&(\psi+S_0)(T_1\sum_{i=0}^nH^{n-i}(H+\psi+S_0)^i+T_2\sum_{i=0}^nH^{n-i}(H+\psi+S_0)^i)& {\rm (by\; (\ref{rel:31}))}\\
 =&S_1((H+\psi+S_0)^{n+1}-H^{n+1})\\
 =&S_1(h_1+S_0)^{n+1},
\end{align*}
where we use the relations, $T_1(\psi+S_0)T_2=0$, $T_1(\psi+S_0)U_2=0$ and $T_2(\psi+S_0)U_1=0$.
The third relation is obtained from multiplying the relation (\ref{rel:5e}) by $(\psi+S_0)$.
\begin{align*}
 0
 =&(\psi+S_0)\cdot (\ref{rel:5e})\\
 =&(\psi+S_0)\cdot\sum_{i=0}^n(H+\psi+S_0)^{n-i}((H+2\psi+2S_0+S_1+S_2)^i-(H+2\psi+2S_0+S_2)^i)\\
 &+(\psi+S_0)\cdot\sum_{i=0}^n(H+\psi)^{n-i}((H+2\psi+S_2)^i-(H+2\psi)^i)\\
 &+2(\psi+S_0)\cdot\sum_{i=0}^nH^{n-i}(H+2\psi+2S_0)^i\\
 =&(\psi+S_0)\cdot\sum_{i=0}^n(H+\psi+S_0)^{n-i}((H+2\psi+2S_0+S_1+S_2)^i-(H+2\psi+2S_0+S_2)^i)\\
 &+(\psi+S_0)\cdot\sum_{i=0}^n(H+\psi+S_0)^{n-i}((H+2\psi+2S_0+S_2)^i-(H+2\psi+2S_0)^i)\\
 &+2(\psi+S_0)\cdot\sum_{i=0}^nH^{n-i}(H+2\psi+2S_0)^i\\
 =&(\psi+S_0)\cdot\sum_{i=0}^n(H+\psi+S_0)^{n-i}(H+2\psi+2S_0+S_1+S_2)^i\\
 &-(\psi+S_0)\cdot\sum_{i=0}^n(H+\psi+S_0)^{n-i}(H+2\psi+2S_0)^i\\
 &+2(\psi+S_0)\cdot\sum_{i=0}^nH^{n-i}(H+2\psi+2S_0)^i\\
 =&(\psi+S_0)\cdot\sum_{i=0}^n(H+\psi+S_0)^{n-i}(H+2\psi+2S_0+S_1+S_2)^i\\
 &+(H+\psi+S_0)^{n+1}-H^{n+1}\\
 =&(\psi+S_0)\cdot\sum_{i=0}^n(H+\psi+S_0)^{n-i}((H+2\psi+2S_0+S_1+S_2)^i+H^i)\\
 =&(g_1-g_0)\sum_{i=0}^ng_1^{n-i}(g_2^i+g_0^i).
\end{align*}
On the second equation, we used $(\psi+S_0)(U_1+T_1)(U_2+T_2)=0$ which comes from (\ref{rel:31}), and $U_1U_2=0$.
On the third equation, we inserted $S_0$ to $(H+\psi)^{n-i}$ since $(H+2\psi+S_2)^i-(H+2\psi)^i$ was divisible by $S_2$, and $S_0S_2=0$.

The last relation is shown by multiplying the second relation by $(g_1-g_2)$ and by applying the first relation to the 
resulting expression.
\begin{align*}
 (g_1-g_0)(g_1-g_2)\sum_{i=0}^n(g_0^i+g_2^i)g_1^{n-i}
 =&(g_1-g_0)(g_1^{n+1}-g_2^{n+1}+(g_1-g_2)\sum_{i=0}^ng_0^ig_1^{n-i})\\
 =&(g_1-g_0)g_1^{n+1}-(g_1-g_0)g_2^{n+1}+(g_1-g_2)g_1^{n+1}\\
 =&-(g_1-g_0)g_2^{n+1}.\quad \Box
\end{align*}
Next, we compute some intersection numbers on $\overline{M}_{0,2}(\bP^n,2)$.
\begin{lemm}\label{lemm:s1}
 \begin{align*}
  &\int S_1^a h_0^bh_1^ch_2^d=0\;{\rm (for}\;\; a+b+c+d=3n+1,\;\;0 < a < n{\rm )}.,\\
  &\int S_1^nh_1h_0^nh_2^n=-1.
 \end{align*}
\end{lemm}
{\bf Proof.} \; Let us discuss the first equation. If $c=0$, then $b+d=3n+1-a>2n+1$. Hence $S_1^a h_0^bh_2^d=0$ by the 
relation $h_0^{n+1}=h_2^{n+1}=0$.
If $c\neq 0$, then thanks to the relation $S_1(h_1-h_0)(h_1-h_2)=0$ in Lemma \ref{lemm:14}, we have only to consider the case of $c=1$. Then $b+d=3n-a>2n$, and $S_1^a h_0^bh_1h_2^d=0$ as above.
As for the second equation, we compute $H^n(H+2\psi+S_2)^n \cdot ((\ref{rel:521})+(\ref{rel:522}))$:
\begin{align*}
 0=&H^n(H+2\psi+S_2)^n \cdot ((\ref{rel:521})+(\ref{rel:522}))& \\
 =&H^n(H+2\psi+S_2)^n (U_1T_2^n+U_2T_1^n+(U_1+U_2)(H+\psi)^n)& {\rm (by\;}H^{n+1}=0{\rm )}\\
 =&h_0^nh_2^n(U_1T_2^n+U_2T_1^n+(h_0-2h_1+h_2)h_1^n)\\
 =&h_0^nh_2^n(U_1T_2^n+U_2T_1^n)-2h_0^nh_1^{n+1}h_2^n.
\end{align*}
It is clear that $U_1T_2^n+U_2T_1^n-S_1^{n-1}P_3=U_1T_2^n+U_2T_1^n-(T_{1}+T_{2})^{n-1}(T_{1}U_{2}+T_{2}U_{1})$ is divisible by $P_1=T_{1}T_{2}$.
Then by using the relation $P_1(h_2-h_0)=0$ of Lemma \ref{lemm:14} and $h_0^{n+1}=h_2^{n+1}=0$,
we obtain $h_0^nh_2^n(U_1T_2^n+U_2T_1^n)=h_{0}^n h_{2}^nS_{1}^{n-1}P_{3}$.
Since $P_3=S_1(h_2-h_1)$, we obtain the following equation:
\begin{equation*}
 -S_1^nh_0^nh_1h_2^n-2h_0^nh_1^{n+1}h_2^n=0.\Box
\end{equation*}
Now, we compute intersection numbers of the type given in (\ref{int num}).
\begin{lemm}\label{lemm:int num}
 Let $\alpha$, $\beta$, $\gamma$, $\delta$ be nonnegative integers which satisfies $\alpha+\beta+\gamma+\delta=3n+1$. If $1\leq \delta\leq n$, then
 \begin{equation*}
  \int g_0^\alpha g_1^\beta g_2^\gamma h_2^\delta=
  \begin{cases}
   \frac{1}{2^{\beta-n}}\Big\{\binom{\beta-n-1}{n-\alpha}-\binom{\beta-n-1}{n-\gamma}\Big\}, (\gamma<n \;{\rm or}\; \gamma=n,\; n-\delta+1>\alpha)\\
   -1, \qquad (\gamma = n,\; n-\delta+1\leq \alpha \leq n)\\
   0, \qquad (0\leq \alpha\leq n,\; n+1\leq \gamma\leq 2n-\delta).
  \end{cases}
 \end{equation*}
 If $\delta=0$, $0\leq \alpha\leq n$, $0\leq \gamma\leq 2n$, then it is zero.
\end{lemm}
{\bf Proof.} \; First, we consider the case of $\delta=0$, $0\leq \alpha\leq n$ and $0\leq \gamma\leq 2n$.
If $n+1\leq \gamma \leq 2n$, then $\alpha+\beta\geq n+1$. Hence it is zero since there is the relations $(g_1-g_0)g_2^{n+1}=0$ and $g_0^{n+1}=0$ of Lemma \ref{lemm:ggg}.
If $0\leq \gamma \leq n+1$, then we have only to consider the case of $g_0^ng_1^{n+1}g_2^n$ since there is the relation $g_1^{n+1}(g_0-2g_1+g_2)=0$ and $g_0^{n+1}=0$.
To prove $g_0^ng_1^{n+1}g_2^n=0$, we should multiply $g_0^ng_2^n$ to the relation $(g_1-g_0)\sum_{i=0}^n(g_0^i+g_2^i)g_1^{n-i}=0$.

Next, we consider the case of $1\leq \delta\leq n-1$, $0\leq \alpha\leq n$, $n+1\leq \gamma\leq 2n-\delta$. In this case, the left hand side is zero since $\alpha+\beta\geq n+1$, and there are relations $(g_1-g_0)g_2^{n+1}=0$ and $g_0^{n+1}=0$.

If $1\leq \delta \leq n$, $\gamma=n$ and $n-\delta+1\leq\alpha\leq n$, then
\begin{align*}
 &h_0^\alpha(h_1+S_0)^\beta(h_2+2S_0+S_1)^n h_2^\delta & \\
 =&h_0^\alpha h_1^\beta(h_2+S_1)^n h_2^\delta & \\
 =&h_0^\alpha h_1^\beta S_1^n h_2^\delta && {\rm (by\; Lemma\; \ref{lemm:s1},\;} n+\delta>n {\rm )}\\
 =&S_1^nh_0^nh_1h_2^n && {\rm (by\; } S_1(h_1-h_0)(h_1-h_2)=0 {\rm )}.
\end{align*}
Therefore, $\int g_0^\alpha g_1^\beta g_2^n h_2^\delta=-1$ by Lemma \ref{lemm:s1}.

Now, we prove the top case.
If $2\leq \delta\leq n$, $n+1-\delta\leq \alpha \leq n$ and $n+1-\delta\leq \gamma \leq n-1$, then the right hand side is zero since $\beta-n-1=2n-\alpha-\gamma-\delta<n-\alpha,n-\gamma$.
On the left hand side, it is
\begin{align*}
 &h_0^\alpha(h_1+S_0)^\beta(h_2+2S_0+S_1)^\gamma h_2^\delta & \\
 =&h_0^\alpha h_1^\beta (h_2+S_1)^\gamma h_2^\delta && {\rm (by}\; \alpha+\delta>n{\rm ,\;\; Lemma\; \ref{rel:mp},\;\ref{lemm:14})}\\
 =&h_0^\alpha h_1^\beta h_2^{\gamma+\delta} && {\rm (by\; Lemma\; \ref{lemm:s1})}\\
 =&0 && {\rm (by}\; \gamma+\delta>n{\rm )}.
\end{align*}

If $\alpha=n$, $\beta=n+1$, $\gamma=n-\delta$ and $1\leq \delta\leq n$, then the right hand side is $\frac12$.
On the left hand side, it is
\begin{align*}
 &h_0^n(h_1+S_0)^{n+1}(h_2+2S_0+S_1)^{n-\delta} h_2^\delta &\\
 =&h_0^nh_1^{n+1}h_2^n&
\end{align*}
Therefore, $\int g_0^ng_1^{n+1}g_2^{n-\delta} h_2^\delta=\frac12$.

If $\alpha=n-\delta$, $\beta=n+1$, $\gamma=n$ and $1\leq\delta\leq n$, then the right hand side is $-\frac12$.
On the left hand side, we consider $0=g_0^{n-\delta}g_2^nh_2^\delta\cdot (g_1-g_0)\sum_{i=0}(g_0^i+g_2^i)g_1^{n-i}$ from Lemma \ref{lemm:ggg}.
\begin{align*}
 0&=g_0^{n-\delta}g_2^nh_2^\delta\cdot (g_1-g_0)\sum_{i=0}^n(g_0^i+g_2^i)g_1^{n-i}&\\
 &=g_0^{n-\delta}g_2^nh_2^\delta\cdot (g_1-g_0)(2g_1^n+\sum_{i=1}^\delta g_0^ig_1^{n-i}) & {\rm (by\; Lemma\; \ref{lemm:ggg})}.
\end{align*}
Here, if $0<i<\delta$, then
\begin{align*}
 &\int g_0^{n-\delta}g_2^nh_2^\delta(g_1-g_0)g_0^ig_1^{n-i}&\\
 =&\int g_0^{n+i-\delta}g_1^{n-i+1}g_2^nh_2^\delta-g_0^{n+i-\delta+1}g_1^{n-i}g_2^nh_2^\delta&\\
 =&(-1)-(-1) & {\rm (by\; the \; above\; case)}&\\
 =&0&.
\end{align*}
Therefore,
\begin{align*}
 0&=\int g_0^{n-\delta}g_2^nh_2^\delta\cdot (g_1-g_0)(2g_1^n+g_0^\delta g_1^{n-\delta}) &\\
 &=\int (2g_0^{n-\delta}g_1^{n+1}g_2^nh_2^\delta-2g_0^{n-\delta}g_1^ng_2^nh_2^\delta+g_0^ng_1^{n-\delta+1}g_2^nh_2^\delta) &\\
 &=\int 2g_0^{n-\delta}g_1^{n+1}g_2^nh_2^\delta -2\cdot(-1)+(-1) &\\
 &=\int 2g_0^{n-\delta}g_1^{n+1}g_2^nh_2^\delta+1.
\end{align*}
Accordingly, $\int g_0^{n-\delta}g_1^{n+1}g_2^n h_2^\delta=-\frac12$.

Finally, we use induction along $\beta\geq n+1$ to prove the remaining part of this lemma.
To execute induction, the following equation is needed: if $1\leq \delta\leq n-1$, $0\leq \alpha \leq n$ and $n+1\leq \gamma\leq 2n+1-\delta$, then
\begin{equation}\label{eq:g2n+1}
  g_0^\alpha g_1^\beta g_2^\gamma h_2^\delta=0.
\end{equation}
It is already proved in the above.

If $\beta=n+1$, then Lemma \ref{lemm:int num} is true as above.
If $\beta>n+1$, by using the relation $g_1^{n+1}(g_0-2g_1+g_2)=0$, we obtain, 
\begin{align*}
 \int g_0^\alpha g_1^\beta g_2^\gamma h_2^\delta
 =&\frac12 \int g_0^{\alpha+1}g_1^{\beta-1}g_2^\gamma h_2^\delta+\frac12 \int g_0^\alpha g_1^{\beta-1}g_2^{\gamma+1} h_2^\delta\\
 =&\frac12 \frac{1}{2^{\beta-n-1}}\Big\{\binom{\beta-n-2}{n-\alpha-1}-\binom{\beta-n-2}{n-\gamma}\Big\}\\
 &+\frac12 \frac{1}{2^{\beta-n-1}}\Big\{\binom{\beta-n-2}{n-\alpha}-\binom{\beta-n-2}{n-\gamma-1}\Big\}\\
 =&\frac{1}{2^{\beta-n}}\Big\{\binom{\beta-n-1}{n-\alpha}-\binom{\beta-n-1}{n-\gamma}\Big\},
\end{align*}
where we use the law of Pascal's triangle in the last line.
The equation \ref{eq:g2n+1} is used in the case of $\gamma=n$ in the above induction.
 $\Box$

In this Lemma, we do not compute the intersection numbers in the cases of $0\leq \delta \leq n$ and $\gamma>2n-\delta$. Although the above result is not complete in this sense, it is sufficient to prove Theorem \ref{theorem:Main2}.

\vspace{1cm}

{\bf Proof of Theorem \ref{theorem:Main2}.}

Let $a$, $b$ be nonnegative integers satisfying $a+b=3n-2k$, and $k$ be a positive integer.
First, we compute intersection number $\int h_0^ah_2^b e^k(g_0,g_1)e^k(g_1,g_2)/(kg_1)$:
\begin{align*}
 \int h_0^ah_2^b \frac{e^k(g_0,g_1)e^k(g_1,g_2)}{kg_1}=&\int h_0^ah_2^b \frac1k \sum_{i=0}^{k-1}\sum_{j=0}^{k-1} \ell_i^k \ell_j^k g_0^{k-i}g_1^{i+j+1}g_2^{k-j}\\
 =&\int \frac1k \sum_{i,j} \ell_i^k \ell_j^k g_0^{a+k-i}g_1^{i+j+1}g_2^{k-j}h_2^b.
\end{align*}
Here we note that it takes a form to which we can apply Lemma \ref{lemm:int num}. Indeed,  if $k$, which is the maximum degree of $g_2$ in this summation, is greater than $2n-b$, then the total degree $a+b+2k+1$ becomes greater than $3n+1$, and the intersection number vanishes. 
We also note that the three cases in Lemma \ref{lemm:int num} are disjoint with each other.
Therefore, we obtain,
\begin{align*}
 &\int h_0^ah_2^b \frac{e^k(g_0,g_1)e^k(g_1,g_2)}{kg_1}\\
 =&\frac1k \sum_{i,j} \ell_i^k \ell_j^k \frac{1}{2^{i+j-n-1}}\Bigg(\binom{i+j-n}{n-a-k+i}-\binom{i+j-n}{n-k+j}\Bigg)\\
 &-\frac1k \sum_{i=k-n+a}^{k-n+a+b-1} \ell_i^k\ell_{k-n}^k\\
 =&\frac1k \sum_{i,j} \ell_i^k \ell_j^k \frac{1}{2^{i+j-n-1}}\Bigg(\binom{i+j-n}{n-a-k+i}-\binom{i+j-n}{n-k+j}\Bigg)\\
 &-\frac1k \ell_{k-n}^k \sum_{i=0}^{b-1} \ell_{a+k-n+i}^k.
\end{align*}
Finally, we compute $\int g_1^ah_2^b e^k(g_0,g_1)e^k(g_1,g_2)/(kg_1)$:
\begin{align*}
 &\int g_1^ah_2^b \frac{e^k(g_0,g_1)e^k(g_1,g_2)}{kg_1}\\
 =&\int \frac1k \sum_{i=0}^{k-1}\sum_{j=0}^{k-1} \ell_i^k \ell_j^k g_0^{k-i}g_1^{a+i+j+1}g_2^{k-j}h_2^b.
\end{align*}
At this stage, we use the fact that the expression in Lemma \ref{lemm:int num}, which contains binomial coefficients, is anti-symmetric
under interchange of $\alpha$ and $\gamma$. Then we obtain,
\begin{align*}
 \int g_1^ah_2^b \frac{e^k(g_0,g_1)e^k(g_1,g_2)}{kg_1}=&\frac1k \ell_{k-n}^k\sum_{i=k-n}^{k-n+b-1}\ell_i^k\\
 =&\frac1k \ell_{k-n}^k\sum_{i=0}^{b-1}\ell_{i+k-n}^k.
\end{align*}
Combining these results, we reach the final expression:
\begin{align*}
 &\int (h_0^a-g_1^a)h_2^b \frac{e^k(g_0,g_1)e^k(g_1,g_2)}{kg_1}\\
 =&\frac1k \sum_{i,j} \ell_i^k \ell_j^k \frac{1}{2^{i+j-n-1}}\Bigg(\binom{i+j-n}{n-a-k+i}-\binom{i+j-n}{n-k+j}\Bigg)\\
 &-\frac1k \ell_{k-n}^k \sum_{i=0}^{b-1} (\ell_{a+k-n+i}^k-\ell_{i+k-n}^k),
\end{align*}
which coincides with the r.h.s. of (\ref{GWinv of deg 2}). $\Box$


\newpage
\begin{thebibliography}{9}
\bibitem{CK0}I. Ciocan-Fontanine, B. Kim, \textit{Wall-crossing in genus zero quasimap theory and mirror maps}, Algebr. Geom. 1 (2014), no. 4, 400-448.
\bibitem{CK} D. A. Cox, S.Kats, \textit{Mirror Symmetry and Algebraic Geometry}, Mathematical Surveys and Monographs Vol.68, American Mathematical Society,1999.
 \bibitem{Cox} D. A. Cox, \textit{The homogeneous coordinate ring of a toric variety, J. Algebraic Geometry}, {\bf 4} (1995), 17-50, alg-geom/9210008.
  \bibitem{JCox} J. A. Cox, \textit{A presentation for the Chow ring $A^*(\bar{M}_{0,2}(P^1,2))$}, Comm. Algebra 35(2007), no. 11, 3391-3414. 
  \bibitem{Ful} W. Fulton, \textit{Introduction to Toric Varieties}, Princeton University Press, Princeton, 1993.
  \bibitem{Ful2} W. Fulton, \textit{Intersection Theory}, Springer-Verlag,1984.
  \bibitem{Givental} A. Givental, \textit{A mirror theorem for toric complete intersections}, Topological field theory, primitive forms and related topics (Kyoto, 1996), 141-175, Progr. Math., 160, Birkh\"{a}user Boston, Boston, MA, 1998.
  \bibitem{Jin1} M. Jinzenji, \textit{Mirror Map as Generating Function of Intersection Numbers: Toric Manifolds with Two K\"{a}hler Forms}, Comm. Math. Phys. 323(2013), no. 2, 747-811.
  \bibitem{Jin2} M. Jinzenji, \textit{Classical Computation of Number of Lines in Projective Hypersurfaces: Origin of Mirror Transformation}, arXiv:1201.5717
  \bibitem{Jin3} M. Jinzenji, \textit{Direct Proof of Mirror Theorem of Projective Hypersurfaces up to degree $3$ Rational Curves}, J. Geom.Phys.61 (2011), no. 8, 1564-1573.
  \bibitem{Kontsevich} M. Kontsevich, \textit{Enumeration of Rational Curves via Torus Actions}, The Moduli space of curves, R.Dijkgraaf, C.Faber, G.van der Geer (Eds.), Progress in Math., v.129, Birkh\"{a}user, 1995, 335-368.
  \bibitem{LLY} B. Lian, K. Liu and S. T. Yau, \textit{Mirror Principle III}, Asian J. Math. {\bf 3} (1999),no.4, 771-800
  \bibitem{Mustata1} A. Musta\c{t}\v{a}, M. Musta\c{t}\v{a}, \textit{Intermediate Moduli Spaces of Stable Maps}, Invent. Math. 167(2007), no. 1,47-90.
  \bibitem{Mustata2} M. Musta\c{t}\v{a}, A. Musta\c{t}\v{a}, \textit{The Chow ring of $\overline{M}_{0,m}({\mathbb P}^n,d)$}, J. Reine Angew. Math. 615 (2008), 93-119.
\end{thebibliography}
\end{document}